\documentclass[a4paper,10pt,
colorlinks,urlcolor=black,linkcolor=black,citecolor=black,
]{article}
\usepackage{amsmath,amssymb}
\usepackage[dvips]{color}
\usepackage{cite}
\usepackage{hangcaption}
\usepackage{amsmath,amssymb}
\usepackage[dvips]{color}
\usepackage{cite}
\usepackage{eepic}
\usepackage{epic}
\usepackage{hangcaption}
\usepackage{amsmath}
\usepackage{amssymb}
\usepackage{amscd}
\usepackage{amsthm}
\usepackage{multicol}
\usepackage[dvipdfmx,%
 bookmarks=true,%
 bookmarksnumbered=true,%
 setpagesize=false,%
 pdfkeywords={TeX; dvipdfmx; hyperref; color;}]{hyperref}
\def\cases{\left\{\begin{array}{ll}}
\def\endcases{\end{array}\right.}
\def\roman{\rm}
\def\bigtimes{\mathop{\mbox{\Large $\times$}}}
\begin{document}
\setcounter{page}{1}
\vskip1.5cm
\begin{center}
{\Large \bf
A quantum linguistic characterization of the reverse relation
between confidence interval and hypothesis testing
}
\vskip0.5cm
{\rm
\large
Shiro Ishikawa
}
\\
\vskip0.2cm
\rm
\it
Department of Mathematics, Faculty of Science and Technology,
Keio University,
\\ 
3-14-1, Hiyoshi, Kouhoku-ku Yokohama, Japan.
E-mail:
ishikawa@math.keio.ac.jp
\end{center}
\par
\rm
\vskip0.3cm
\par
\noindent
{\bf Abstract}
\normalsize
\vskip0.5cm
\par
\noindent
Although there are many ideas for the formulations of statistical hypothesis testing,
we consider that the likelihood ratio test is the most reasonable and orthodox. However, it is not handy, and thus, it is not usual in elementary books. That is,
the statistical hypothesis testing written in elementary books is different 
from the likelihood ratio test.
Thus, from the theoretical point of view, we have the following question:
\begin{itemize}
\item{}
What is the statistical hypothesis testing written in elementary books?
\end{itemize}
For example, we consider that even the difference between "one sided test" and "two sided test" is not clear yet. In this paper, we give an answer to this question.
That is,
we propose a new formulation of statistical hypothesis testing,
which is contrary to the confidence interval methods. In other words, they are two sides of the same coin. This will be done in quantum language (or, measurement theory), which is characterized as the linguistic turn of the Copenhagen interpretation of quantum mechanics, and also, a kind of system theory such that it is applicable to both classical and quantum systems.
Since quantum language is suited for theoretical arguments, we believe that our results are essentially final as a general theory.

%
%
%

\vskip1.0cm

\par

\def\Cal{\cal}
\def\bigstimes{\text{\large $\: \boxtimes \,$}}

\par
\noindent

\vskip0.2cm
\par
\noindent
\par
\noindent
\section{
Introduction
}
%

\rm
\par
\par
\noindent

%

\par
\noindent
\subsection{
Quantum language
(Axioms
and
Interpretation)
}
As mentioned in the above abstract, our purpose is to answer the following question:
\begin{itemize}
\item[(A)]
What is the statistical hypothesis testing written in elementary books?
\end{itemize}
This will be answered in terms of quantum language.
\par
According to ref.\cite{Ishi8},
we shall mention the overview of quantum language
(or, measurement theory, in short, MT).
\par
\par
\rm
Quantum language is characterized as the linguistic turn of the Copenhagen interpretation of quantum mechanics({\it cf.} ref.(\cite{Ishi5},
{{{}}}{\cite{Neum}}).
Quantum language (or, measurement theory ) has two simple rules
(i.e. Axiom 1(concerning measurement) and Axiom 2(concerning causal relation))
and the linguistic interpretation (= how to use the Axioms 1 and 2). 
That is,
\begin{align}
\underset{\mbox{(=MT(measurement theory))}}{\fbox{Quantum language}}
=
\underset{\mbox{(measurement)}}{\fbox{Axiom 1}}
+
\underset{\mbox{(causality)}}{\fbox{Axiom 2}}
+
\underset{\mbox{(how to use Axioms)}}{\fbox{linguistic interpretation}}
\label{eq1}
\end{align}
({\it cf.} refs.
{{{}}}{\cite{Ishi2}-\cite{Ishi9}}).
\par
This theory is formulated in a certain $C^*$-algebra ${\cal A}$({\it cf.} ref.
{{{}}}{\cite{Saka}}), and is classified as follows:
\begin{itemize}
\item[(B)]
$
\quad
\underset{\text{\scriptsize }}{\text{MT}}
$
$\left\{\begin{array}{ll}
\text{quantum MT$\quad$(when ${\cal A}$ is non-commutative)}
\\
\\
\text{classical MT
$\quad$
(when ${\cal A}$ is commutative, i.e., ${\cal A}=C_0(\Omega)$)}
\end{array}\right.
$
\end{itemize}
where $C_0(\Omega)$
is
the $C^*$-algebra composed of all continuous 
complex-valued functions vanishing at infinity
on a locally compact Hausdorff space $\Omega$.

Since our concern in this paper is concentrated to 
the usual statistical hypothesis test
methods in statistics,
we devote ourselves to the commutative $C^*$-algebra $C_0(\Omega)$,
which is quite elementary.
Therefore, we believe that all statisticians
can understand our assertion
(i.e.,
a new viewpoint of the confidence interval
methods
).

Let $\Omega$ is a locally compact Hausdorff space, which is also called
a state space. And thus, an element $\omega (\in \Omega )$ is said to be a state.
Let $C(\Omega)$ be the $C^*$-algebra composed of all bounded continuous 
complex-valued functions on a locally compact Hausdorff space $\Omega$.
The norm $\| \cdot \|_{C(\Omega )}$ is usual, i.e.,
$\| f \|_{C(\Omega )} = \sup_{\omega \in \Omega } |f(\omega )|$
$(\forall f \in C(\Omega ))$.
\par
\rm
Motivated by Davies' idea ({\it cf.} ref.
{{{}}}{\cite{Davi}}) in quantum mechanics,
an observable ${\mathsf O}=(X, {\mathcal F}, F)$ in $C_0(\Omega )$
(or, precisely, in $C(\Omega )$) is defined as follows:
\begin{itemize}
\item[(C$_1$)]
$X$ is a topological space. 
${\mathcal F} ( \subseteq 2^X$(i.e., the power set of $X$) is a field,
that is, it satisfies the following conditions (i)--(iii):
(i):
$\emptyset \in {\cal F}$, 
(ii):$\Xi \in {\mathcal F} \Longrightarrow X\setminus \Xi  \in 
{\mathcal F}$,
(iii):
$\Xi_1, \Xi_2,\ldots, \Xi_n \in {\mathcal F} \Longrightarrow \cup_{k=1}^n \Xi_k \in {\mathcal F}$.
\item[(C$_2$)]
The map $F: {\cal F} \to C(\Omega )$ satisfies that 
\begin{align}
0 \le [F(\Xi )](\omega ) \le 1, \quad [F(X )](\omega )=1
\qquad
(\forall \omega \in \Omega )
\nonumber \end{align}
and moreover, 
if
\begin{align}
\Xi_1, \Xi_2,\ldots, \Xi_n, \ldots \in {\mathcal F},
\quad
\Xi_m \cap \Xi_n = \emptyset \quad( m \not= n ),
\quad
\Xi = \cup_{k=1}^\infty \Xi_k \in {\mathcal F},
\nonumber \end{align} 
then, it holds
\begin{align}
[F(\Xi)](\omega) = \lim_{n \to \infty } \sum_{k=1}^n [F(\Xi_k )](\omega )
\quad
(\forall \omega \in \Omega )
\nonumber \end{align}
\end{itemize}
Note that Hopf extension theorem
({\it cf.}
ref.
{{{}}}{\cite{Yosi}})
guarantees that
$(X, {\cal F}, [F(\cdot)](\omega))$
is regarded as the mathematical probability space.
\par
\noindent
\bf
Example 1
\rm
[Normal observable].
Let ${\mathbb R}$ be the set of the real numbers.
Consider the state space
$\Omega = {\mathbb R} \times {\mathbb R}_+$,
where 
${\mathbb R}_+=\{ \sigma \in {\mathbb R} | \sigma > 0 \}$.
Define the normal observable
${\mathsf O}_N = ({\mathbb R}, {\mathcal B}_{\mathbb R}, {{{N}}} )$ in $C_0({\mathbb R} \times {\mathbb R}_+)$ such that
\begin{align}
&
[{{{N}}}({\Xi})] ({} {}{\omega} {}) 
=
\frac{1}{{\sqrt{2 \pi }\sigma{}}}
\int_{{\Xi}} \exp[{}- \frac{({}{}{x} - {}{\mu}  {})^2 }{2 \sigma^2}    {}] d {}{x}
\label{eq2}
\\
&
\quad
({}\forall  {\Xi} \in {\cal B}_{{\mathbb R}{}}\mbox{(=Borel field in ${\mathbb R}$))},
\quad
\forall   {}{\omega} =(\mu, \sigma)   \in \Omega = {\mathbb R}{}\times {\mathbb R}_+).
\nonumber
\end{align}
In this paper, we devote ourselves to the normal observable.

\vskip0.5cm

%
%

\par
Now we shall briefly explain "quantum language (\ref{eq1})" in classical systems as follows:
A measurement of an observable
${\mathsf O}=(X, {\mathcal F}, F)$
for a system with a state $\omega (\in \Omega )$
is denoted by
${\mathsf M}_{C_0(\Omega)} ({\mathsf O}, S_{[\omega]})$.
By the measurement, a measured value $x (\in X)$ is obtained as follows:
\par
\noindent
\bf
Axiom 1
\rm
(Measurement)
\begin{itemize}
\item{}
\sl
The probability that a measured value $x$
$( \in X)$ obtained by the measurement 
${\mathsf{M}}_{{{C_0(\Omega)}}} ({\mathsf{O}}$
${ \equiv} (X, {\cal F}, F),$
{}{$ S_{[\omega_0]})$}
belongs to a set 
$\Xi (\in {\cal F})$ is given by
$
[F(\Xi) ](\omega_0 )
$.
\end{itemize}
\rm
\par
\noindent
\par
\noindent
\bf
Axiom 2
\rm
(Causality)
\begin{itemize}
\item{}
\sl
The causality is represented by a Markov operator
$\Phi_{21} : C_0(\Omega_2 ) \to C_0(\Omega_1 )$.
Particularly, the deterministic causality
is represented by a continuous map
$\pi_{12} : \Omega_1 \to \Omega_2$
\end{itemize}
\par
\noindent
\bf
Interpretation
\rm
(Linguistic interpretation).
Although there are several linguistic rules in quantum language, the following is the most important:
\begin{itemize}
\item{}
\sl
Only one measurement is permitted.
\end{itemize}
\rm
In order to read this paper,
it suffices to understand the above three.
For the further arguments, see refs.
{{{}}}{\cite{Ishi2}-\cite{Ishi9}}.
\vskip0.5cm

Consider measurements ${\mathsf{M}}_{{{C_0(\Omega)}}} ({\mathsf{O}_k}$
${ \equiv} (X_k, {\cal F}_k, F_k),$
{}{$ S_{[\omega_0]})$}, $(k=1,2, \ldots, n )$. However, the linguistic interpretation says that
only one measurement is permitted.
Thus we must consider a simultaneous measurement or a parallel measurement.

\par
\noindent
\bf
Definition 1
\rm
[(i):Simultaneous observable].
\sl
Let
${\mathsf{O}_k}$
${ \equiv} (X_k, {\cal F}_k, F_k)$
$(k=1,2, \ldots, n)$
be an observable in $C_0(\Omega)$.
The simultaneous observable
$\bigtimes_{k=1}^{n}{\mathsf{O}_k}$
${ \equiv} (\bigtimes_{k=1}^n X_k, \bigstimes_{k=1}^n {\cal F}_k, \widehat{F}
(\equiv \bigtimes_{k=1}^n F_k
))$
in $C_0(\Omega)$ is defined by
\begin{align}
&
[\widehat{F}(\Xi_1 \times \cdots \times \Xi_n )](\omega)
(\equiv
[(\bigtimes_{k=1}^n F_k)(\Xi_1 \times \cdots \times \Xi_n )](\omega)
)
=
\bigtimes_{k=1}^n [F_k(\Xi_k)](\omega)
\label{eq3}
\\
&
\quad
\qquad
\qquad
(\forall \Xi_k \in {\mathcal F}_k \;\;(k=1, \ldots, n), \forall \omega \in \Omega )
\nonumber
\end{align}
Here,
$ \boxtimes_{k=1}^n {\cal F}_k$
is the smallest field including
the family
$\{
{\text{\large $\times$}}_{k=1}^n \Xi_k
$
$:$
$\Xi_k \in {\cal F}_k \; k=1,2,\ldots, n \}$.
If 
${\mathsf{O}}$
${ \equiv} (X, {\cal F}, F)$
is equal to
${\mathsf{O}_k}$
${ \equiv} (X_k, {\cal F}_k, F_k)$
$(k=1,2, \ldots, n)$,
then
the simultaneous observable
$\bigtimes_{k=1}^{n}{\mathsf{O}_k}$
${ \equiv} (\bigtimes_{k=1}^n X_k, \bigstimes_{k=1}^n {\cal F}_k, \widehat{F}
(\equiv \bigtimes_{k=1}^n F_k
))$
is denoted by
${\mathsf{O}^n}$
${ \equiv} (X^n, {\cal F}^n, F^n)$.
\par
\noindent
\rm
[(ii):Parallel observable].
\sl
Let
${\mathsf{O}_k}$
${ \equiv} (X_k, {\cal F}_k, F_k)$
be an observable in $C_0(\Omega_k)$,
$(k=1,2, \ldots, n)$.
The parallel observable
$\bigotimes_{k=1}^{n}{\mathsf{O}_k}$
${ \equiv} (\bigtimes_{k=1}^n X_k, \bigstimes_{k=1}^n {\cal F}_k, \widetilde{F}
(\equiv \bigotimes_{k=1}^n F_k
))$
in $C_0(\bigtimes_{k=1}^n \Omega_k)$ is defined by
\begin{align}
&
[\widetilde{F}(\Xi_1 \times \cdots \times \Xi_n )](\omega_1,\omega_2, \ldots,\omega_n)
(\equiv
[(\bigotimes_{k=1}^n F_k)(\Xi_1 \times \cdots \times \Xi_n )](\omega_1,\omega_2, \ldots,\omega_n)
)
=
\bigtimes_{k=1}^n [F_k(\Xi_k)](\omega_k)
\label{eq4}
\\
&
\quad
\qquad
\qquad
\quad
\qquad
\qquad
(\forall \Xi_k \in {\mathcal F}_k , \forall \omega_k \in \Omega_k,\;\;(k=1, \ldots, n))
\nonumber
\end{align}

\par
\noindent
\bf
Definition 2
\rm
[Image observable].
\sl
Let
${\mathsf{O}}$
${ \equiv} (X, {\cal F}, F)$
be observables in $C_0(\Omega)$.
The observable $f({\mathsf{O}})$ 
$({ \equiv} (Y, {\cal G}, G (\equiv F \circ f^{-1}))$
in $C_0(\Omega)$
is called the image observable of
${\mathsf{O}}$
by a map $f:X \to Y$,
if it holds that
\begin{align}
G( \Gamma ) = F( f^{-1}(\Gamma))
\qquad
(
\forall \Gamma \in {\mathcal G}
)
\label{eq5}
\end{align}

\par
\noindent
\bf
Example 2
\rm
[Simultaneous normal observable].
Let
${\mathsf O}_N = ({\mathbb R}, {\mathcal B}_{\mathbb R}, {{{N}}} )$ be the normal observable
in $C_0({\mathbb R} \times {\mathbb R}_+)$ in Example 1.
Let $n$ be a natural number.
Then, we get the simultaneous normal observable
${\mathsf O}_N^n = ({\mathbb R}^n, {\mathcal B}_{\mathbb R}^n, {{{N}}^n} )$ 
in $C_0({\mathbb R} \times {\mathbb R}_+)$.
That is,
\par
\noindent
\begin{align}
&
[{{{N}}}^n
(\bigtimes_{k=1}^n \Xi_k)]
({}\omega{})
=
\bigtimes_{k=1}^n
[{{{N}}}(\Xi_k)](\omega)
\nonumber
\\
=
&
\frac{1}{({{\sqrt{2 \pi }\sigma{}}})^n}
\underset{{\bigtimes_{k=1}^n \Xi_k }}{\int \cdots \int}
\exp[{}- \frac{\sum_{k=1}^n ({}{}{x_k} - {}{\mu}  {})^2 
}
{2 \sigma^2}    {}] d {}{x_1} d {}{x_2}\cdots dx_n
\label{eq6}
\\
&
\qquad 
({}\forall  \Xi_k \in {\cal B}_{{\mathbb R}{}}^{}
({}k=1,2,\ldots, n),
\quad
\forall   {}{\omega}=(\mu, \sigma )    \in \Omega = {\mathbb R}\times {\mathbb R}_+{}).
\nonumber
\end{align}

Consider the maps
$\overline{\mu}: {\mathbb R}^n \to {\mathbb R}$, 
${\overline{SS}}: {\mathbb R}^n \to {\mathbb R}$
and
${\overline{\sigma}}: {\mathbb R}^n \to {\mathbb R}$
such that
\begin{align}
&
\overline{\mu}
(x) =
\overline{\mu}
(x_1,x_2,\ldots , x_n ) =
\frac{x_1 + x_2 + \cdots + x_n}{n}
\quad( \forall x=(x_1,x_2,\ldots , x_n ) \in {\mathbb R}^n )
\label{eq7}
\\
&
{{\overline{SS}}}
(x) =
{{\overline{SS}}}
(x_1,x_2,\ldots , x_n ) =
{\sum_{k=1}^n ( x_k - 
\overline{\mu}
(x))^2}
\quad( \forall x=(x_1,x_2,\ldots , x_n ) \in {\mathbb R}^n)
\nonumber
\\
&
{{\overline{\sigma}}}
(x) =
{{\overline{\sigma}}}
(x_1,x_2,\ldots , x_n ) =
\sqrt{
\frac
{\sum_{k=1}^n ( x_k - 
\overline{\mu}
(x))^2}
{n}
}
\quad( \forall x=(x_1,x_2,\ldots , x_n ) \in {\mathbb R}^n)
\label{eq8}
\end{align}
Thus, we have two
image observables
$\overline{\mu}({\mathsf O}_N^n) $
$= ({\mathbb R}, {\mathcal B}_{\mathbb R}, {{{N}}^n} \circ \overline{\mu}^{-1} )$
and
${{\overline{SS}}}({\mathsf O}_N^n) $
$= ({\mathbb R}_+, {\mathcal B}_{{\mathbb R}_+}, {{{N}}^n} \circ {{\overline{SS}}}^{-1} )$
in $C_0({\mathbb R} \times {\mathbb R}_+)$.

It is easy to see that
\begin{align}
&
[({{{N}}^n} \circ \overline{\mu}^{-1})(\Xi_1)](\omega)
=
\frac{1}{({{\sqrt{2 \pi }\sigma{}}})^n}
\underset{
\{ x \in {\mathbb R}^n \;:\; {\overline{\mu}}(x) \in \Xi_1 \}}
{\int \cdots \int}
\exp[{}- \frac{\sum_{k=1}^n ({}{}{x_k} - {}{\mu}  {})^2 
}
{2 \sigma^2}    {}] d {}{x_1} d {}{x_2}\cdots dx_n
\nonumber
\\
=
&
\frac{\sqrt{n}}{{\sqrt{2 \pi }\sigma{}}}
\int_{{\Xi_1}} \exp[{}- \frac{n({}{}{x} - {}{\mu}  {})^2 }{2 \sigma^2}    {}] d {}{x}
\label{eq9}
\intertext{and}
&
[({{{N}}^n} \circ {{{\overline{SS}}}}^{-1})(\Xi_2)](\omega)
=
\frac{1}{({{\sqrt{2 \pi }\sigma{}}})^n}
\underset{
\{ x \in {\mathbb R}^n \;:\; {\overline{SS}}(x) \in \Xi_2 \}}
{\int \cdots \int}
\exp[{}- \frac{\sum_{k=1}^n ({}{}{x_k} - {}{\mu}  {})^2 
}
{2 \sigma^2}    {}] d {}{x_1} d {}{x_2}\cdots dx_n
\nonumber
\\
=
&
\int_{\Xi_2 / \sigma^2} p^{{\chi}^2}_{n-1}({ x} ) {dx}
\label{eq10}
\\
&
\quad
({}\forall  {\Xi_1} \in {\cal B}_{{\mathbb R}{}},
\;\;
\forall \Xi_2 \in {\cal B}_{{\mathbb R}_+{}},
\quad
\forall   {}{\omega} =(\mu, \sigma)   \in \Omega \equiv {\mathbb R}{}\times {\mathbb R}_+).
\nonumber
\end{align}
Here, $p^{{\chi}^2}_{n-1}({ x} )$ is the chi-squared distribution with $n-1$ degrees of freedom. That is,
\begin{align}
p^{{\chi}^2}_{n-1}({ x} )
=
\frac{x^{(n-1)/2-1}e^{-x/2}}{2^{(n-1)/2} \Gamma ((n-1)/2)}
\quad ( x > 0)
\label{eq11}
\end{align}
where $\Gamma$ is the gamma function.

%
%

\subsection{
Fisher's maximum likelihood method}

\rm
\par
\noindent 
\par
It is usual to consider that
we do not know the pure state
$\omega_0$
$(
\in
\Omega
)$
when
we take a measurement
${\mathsf{M}}_{{{C_0(\Omega)}}} ({\mathsf{O}}, S_{[\omega_0]})$.
That is because
we usually take a measurement ${\mathsf{M}}_{{{C_0(\Omega)}}} ({\mathsf{O}},
S_{[\omega_0]})$
in order to know the state $\omega_0$.
Thus,
when we want to emphasize that
we do not know the state $\omega_0$,
${\mathsf{M}}_{{{C_0(\Omega)}}} ({\mathsf{O}}, S_{[\omega_0]})$
is denoted by
${\mathsf{M}}_{{{C_0(\Omega)}}} ({\mathsf{O}}, S_{[\ast]})$.
Also,
if
we know (or, postulate) that a state $\omega_0$
belongs to
a certain suitable set $K$
$(\subseteq 
\Omega )$,
the
${\mathsf{M}}_{{{C_0(\Omega)}}} ({\mathsf{O}}, S_{[\omega_0]})$
is denoted by
\begin{align}
{\mathsf{M}}_{{{C_0(\Omega)}}} ({\mathsf{O}}, S_{[\ast]}
(K) ).
\label{eq100}
\end{align}

\vskip0.5cm
\par
\noindent
{\bf Theorem 1}
\rm
[Fisher's maximum likelihood method ({\it cf.} refs.
{{{}}}{\cite{Ishi3},\cite{Ishi4})].
\sl
Consider a measurement
${\mathsf M}_{{C_0(\Omega)}}
(
{\mathsf O}=(X , {\cal F} , F )
,$
$ S_{[*]}(K))$.
Assume that
we know that the measured value $x \;(\in X )$
obtained by a measurement
${\mathsf M}_{{C_0(\Omega)}}
(
{\mathsf O}=(X , {\cal F} , F )
,$
$ S_{[*]}(K))$
belongs to
$\Xi (\in {\cal F})$.
Then,
there is a reason to infer that the unknown state
$[\ast ]$ is equal to $\omega_0 (\in K )$ such that
\begin{align}
\min_{\omega_1 \in K}
\frac{[F(\Xi)](\omega_0)}
{
[F(\Xi)](\omega_1)
}
\Big(
=
\frac{[F(\Xi)](\omega_0)}
{
\max_{\omega_1 \in K} [F(\Xi)](\omega_1)
}
\Big)
=
1
\label{eq12}
\end{align}
if the righthand side of this formula exists.
Also, if $\Xi=\{x\}$, it suffices to calculate the $\omega_0 (\in K)$ 
such that
$$
L(x, \omega_0)
=1
$$
where the likelihood function $L(x, \omega ) (\equiv L_x (\omega ))$ is defined by
\begin{align}
L(x, \omega )
=
\inf_{\omega_1 \in K }
\Big[\lim_{\Xi \supseteq \{x \}, \; [F(\Xi)](\omega_1) \not= 0,\; \Xi \to \{x \}} \frac{[F(\Xi)](\omega)}{
[F(\Xi)](\omega_1)}
\Big]
\label{eq13}
\end{align}
\par
\noindent
\rm
\bf
Example 3
\rm
[Fisher's maximum likelihood method].
Consider the simultaneous normal observable
${\mathsf O}_N^n = ({\mathbb R}^n, {\mathcal B}_{\mathbb R}^n, {{{N}}^n} )$
in $C_0({\mathbb R} \times {\mathbb R}_+)$ in the formula (\ref{eq6}).
Thus, we have the simultaneous measurement
${\mathsf M}_{C_0({\mathbb R} \times {\mathbb R}_+ )} ({\mathsf O}_N^n = ({\mathbb R}^n, {\mathcal B}_{\mathbb R}^n, {{{N}}^n} )$,
$S_{[\ast]}(K))$ 
in $C_0({\mathbb R} \times {\mathbb R}_+)$.
Assume that a measured value $x=(x_1, x_2, \ldots, x_n ) (\in
{\mathbb R}^n )$ is obtained by the measurement.
Since the likelihood function
$L_x(\mu, \sigma)(=L(x, (\mu,\sigma)) $
is defined by
\begin{align}
L_x(\mu, \sigma)
&
=
\frac{1}{({{\sqrt{2 \pi }\sigma{}}})^n}
\exp[{}- \frac{\sum_{k=1}^n ({}{}{x_k} - {}{\mu}  {})^2 
}
{2 \sigma^2}    {}] 
\nonumber
\intertext{or, in the sense of (\ref{eq13}),}
L_x(\mu, \sigma)
&
=
\frac{
\frac{1}{({{\sqrt{2 \pi }\sigma{}}})^n}
\exp[{}- \frac{\sum_{k=1}^n ({}{}{x_k} - {}{\mu}  {})^2 
}
{2 \sigma^2}    {}] }
{
\frac{1}{({{\sqrt{2 \pi }\overline{\sigma}(x){}}})^n}
\exp[{}- \frac{\sum_{k=1}^n ({}{}{x_k} - {}{\overline{\mu}(x)}  {})^2 
}
{2 \overline{\sigma}(x)^2}    {}] 
}
\label{eq14}
\\
&
({}\forall x = (x_1, x_2, \ldots , x_n ) \in {\mathbb R}^n,
\quad
\forall   {}{\omega}=(\mu, \sigma )    \in \Omega = {\mathbb R}\times {\mathbb R}_+{}).
\nonumber
\end{align}
it suffices to calculate the following equations:
\begin{align}
\frac{\partial L_x(\mu, \sigma)}{\partial \mu}=0,
\quad
\frac{\partial L_x(\mu, \sigma)}{\partial \sigma}=0
\label{eq15}
\end{align}
\par
\noindent
For example, assume that $K={\mathbb R} \times {\mathbb R}_+ $.
Solving the equation (\ref{eq15}), we can infer, by Theorem 1 (Fisher's maximum likelihood method), that
$[\ast]=(\mu, \sigma)$ 
$(\in 
{\mathbb R} \times {\mathbb R}_+
)$
such that
\begin{align}
\mu=\overline{\mu}(x) =\frac{x_1 + x_2+ \ldots + x_n }{n},
\quad
\sigma=
\overline{\sigma}(x)
=
\sqrt{\frac{\sum_{k=1}^n (x_k - \overline{\mu}(x))^2}{n}}
=
\sqrt{\frac{n-1}{n}}\overline{\sigma}'(x)
\label{eq16}
\end{align}
\par
\noindent


\par
\par
\noindent
\par
\rm

\noindent
\subsection{
The orthodox characterization of statistical hypothesis testing
(the likelihood ratio test)
}

Our purpose of this paper is to
propose a kind of statistical hypothesis test which is characterized as "the
reverse confidence reverse"
in the following Section 2.
However, before it, we mention the standard statistical hypothesis test
(i.e.,
the likelihood ration test)
as follows.

\rm
Consider a measurement
%
${\mathsf M}_{C_0(\Omega)}({\mathsf O}\equiv(X, {\cal F}, F{}), S_{[*]}
)$
formulated in
${C_0(\Omega)}$.
Here, we assume that
$(X,
\tau{{}_X})$ is a topological space,
where
$\tau{{}_X}$
is the set of all open sets.
And assume that
$\overline{\cal F}={\cal B}_X$;
the Borel field,
i,e.,
the smallest $\sigma$-field that contains all
open sets in $X$.
Note that we can assume, without loss of generality,
that
$
F({\Xi})
\not=
0
$
for any open set $\Xi (\in \tau{{}_X} )$
such that $\Xi \not= \emptyset$.
That is because,
if
$
F({\Xi})
=
0
$,
it suffices to
redefine
$X$ by $X\setminus \Xi$.
Let $\Theta$ be a locally compact space
with the Borel field ${\cal B}_\Theta$.
Let
$\pi:\Omega \to \Theta$
be a continuous map,
which is a kind of causal relation (in Axiom 2), and called
\it
{\lq\lq}quantity{\rq\rq},$\;\;$
\rm
and let
$E:X \to \Theta$
be a continuous (or more generally, measurable) map,
which is called
\it
{\lq\lq}estimator{\rq\rq}$.\; \;$
\rm

Assume the following hypothesis called {\lq\lq}{\it null hypothesis}":
\begin{itemize}
\item[{}{(D)}]
$\pi( \ast )$
(where
$[\ast]$
is
the unknown state
in
${\mathsf M}_{C_0(\Omega)}({\mathsf O}, S_{[*]}
)$
)
belongs to
a set ${{H_N}}$
$({}\subseteq 
\Theta
)$.
\end{itemize}
In short, the set $H_N$ is also called {\lq\lq}{\it null hypothesis}".

%

\par
\noindent
In order to deny this hypothesis {}{(D)},
we define
the rejection region $
{\widehat R}^\alpha_{{{H_N}}}
$
($\in {\cal B}_\Theta$) as follows.

\begin{itemize}
\item[{}{\roman{(E)}}]
For sufficiently small
{\it significance level}
$\alpha$
(
$0 < \alpha \ll 1$
,
e.g.,
$\alpha=0.05$
),
define the
{\it
rejection region
}
${\widehat R}^\alpha_{{{H_N}}} \in {\cal B}_\Theta$
such that
\begin{itemize}
\item[(E$_1$)]
$
\sup_{
\omega \in \pi^{-1}( \{ \theta \} )
}
$
$
[F(E^{-1}({\widehat R}^\alpha_{{{H_N}}}))](\omega)
{{{}}}
\le
\alpha
\quad
(\forall 
\theta \in H_N ( \subseteq  \Theta )
)
$
\item[(E$_2$)]
If 
${\widehat R}^{\alpha,1}_{{{H_N}}}
(\in {\cal B}_\Theta )
$
and
${\widehat R}^{\alpha,2}_{{{H_N}}}
(\in {\cal B}_\Theta )$
satisfy
(E$_1$)
and
${\widehat R}^{\alpha,1}_{{{H_N}}} \subseteq 
{\widehat R}^{\alpha,2}_{{{H_N}}}
$,
then,
choose
${\widehat R}^{\alpha,2}_{{{H_N}}}
$.
\end{itemize}
\end{itemize}

\par
\noindent
%

\begin{center}
\unitlength=0.20mm
\begin{picture}(400,130)
\put(-30,0){{
\put(27,18){0}
\put(27,108){1}
\put(350,18){$
\Theta
$}
\dottedline{3}(40,110)(340,110)
\put(40,20){\line(0,1){100}}
\thicklines
\put(40,20){\line(1,0){300}}
\thicklines
\spline(40,110)(60,109)(80,105)(110,100)
(150,30)(180,23)(200,27)(250,25)
(270,35)(280,70)(300,100)(340,110)
\dottedline{5}(40,40)(340,40)
\put(27,40){$\alpha$}
\multiput(150,18)(3,0){40}{\line(0,2){4}}
\put(200,2){${{{H_N}} } $}
\put(134,130){
\scriptsize
$
\footnotesize{
{{{}}}
\underset{
\omega \in \pi^{-1}( \{ \theta \} )
}{\sup}
[
F(
E^{-1}(
{\widehat R}^\alpha_{{{H_N}}}
))](\omega)
{{{}}}
}
$}
}}
\put(40,-20){\bf Figure 1.
\rm
Null Hypothesis 
${{{H_N}}}$
}
\end{picture}
\end{center}
\vskip0.5cm
Then, Axiom 1 says that
\begin{itemize}
\item[(F)]
if
$\pi(\ast) \in {{H_N}}$,
the following (F$_1$)
(or, equivalently, (F$_2$)
)
holds:
\begin{itemize}
\item[(F$_1$)]
the probability that
a measured value
obtained by
${\mathsf M}_{{C_0(\Omega)}}({\mathsf O} $
$\equiv(X, {\cal F}, F{}) , S_{[\ast]})$
belong to
$E^{-1}({\widehat R}^\alpha_{{{H_N}}})$
is less than or equal to
$\alpha$.
\item[(F$_2$)]
the probability that
a measured value
obtained by
${\mathsf M}_{{C_0(\Omega)}}(E{\mathsf O} $
$\equiv(\Theta, {\cal B}_\Theta, F\circ E^{-1}{}) , S_{[\ast]})$
belong to
${\widehat R}^\alpha_{{{H_N}}}$
is less than or equal to
$\alpha$.
\end{itemize}
Therefore,
if
$\pi(\ast) \in {{H_N}}$,
and
if $\alpha$ is sufficiently small,
then
there is a reason to deny the hypothesis
{}{(D)}.
\end{itemize}


\par
\noindent
It is clear
that
the rejection region ${\widehat R}^\alpha_{{{H_N}}}$
is not uniquely determined
in general.
%
%
Thus, we have the following problem:
\begin{itemize}
\item[(G)]
Find the most proper rejection region
${\widehat R}^\alpha_{{{H_N}}}$.
\end{itemize}
In what follows, we shall answer this (G)
as "the likelihood ratio test".

Let
$
E({\mathsf O} )$
$\equiv(\Theta, {\cal B}_\Theta, F\circ E^{-1}{}) $
be
the image observable of
the
${\mathsf O} $
$\equiv$
$(X , {\cal F} , {F})$
in
a
commutative $C^*$-algebra
${C_0(\Omega)}$.
Define
the likelihood function
${L}: \Theta \times 
\Omega
\to
[0,1]$
of
the image observable
$
E( {\mathsf O} )$
by
(\ref{eq13}).
\rm
Let
${{H_N}}$
be 
as in (D).
Here define the
function
$\Lambda_{{{H_N}}}{}: \Theta \to [0,1]$
such that:
\begin{align}
\Lambda_{{{H_N}}}({}\theta)
=
\sup_{ \omega \in \Omega \mbox{ such that } \pi (\omega ) \in {{H_N}}}
L(\theta, \omega )
\quad
({}\forall \theta \in \Theta{}).
\label{eq17}
\end{align}
Also, for any
$\epsilon \; ({}0 < \epsilon \le 1{})$,
define
${{{R}}}_{{{H_N}}}^\epsilon $
$({}\in {\cal B}_\Theta{})$
such that
\begin{align}
{{{R}}}_{{{H_N}}}^\epsilon
=
\{ \theta \in \Theta \; | \;
\Lambda_{{{H_N}}} ({}\theta{}) \le \epsilon \}.
\label{eq18}
\end{align}
\par
\noindent
\begin{center}
\unitlength=0.20mm
\begin{picture}(400,130)
\put(27,18){0}
\put(27,50){$\epsilon$}
\put(27,108){1}
\put(350,18){$\Theta$}
\dottedline{3}(40,110)(340,110)
\put(40,20){\line(0,1){100}}
\thicklines
\put(40,20){\line(1,0){300}}
\thicklines
\spline(40,70)(60,95)(75,108)(80,110)(85,110)(95,110)(100,110)
(150,100)(200,60)(250,40)
(270,35)(280,30)(300,25)(340,20)
\dottedline{5}(40,50)(340,50)
\dottedline{5}(225,50)(225,20)
\multiput(225,18)(3,0){40}{\line(0,2){4}}
\put(280,2){${{{R}}}_{{{H_N}}}^\epsilon $}
\put(190,80){$\Lambda_{{{H_N}}} (\theta) $}
\put(80,-20){\bf Figure 2.
\rm
${R}_{{{H_N}}}^\epsilon$
}
\end{picture}
\end{center}
\vskip0.5cm

\par
\noindent
Consider a positive number
$\alpha$
(called
\it
a significance level
\rm
)
such that
$0 < \alpha \ll 1$
({}e.g.
$\alpha = 0.05 $
).
\rm
Thus we can define
$\epsilon(\alpha)$
such that:
\begin{align}
\epsilon(\alpha)
=
\sup
\{
\epsilon \; |
\;
\sup_{ \omega \in \Omega \mbox{ such that } \pi (\omega ) \in {{H_N}}}
{{{}}} 
[
F(
E^{-1}(
{{{R}}}_{{{H_N}}}^\epsilon{}
)
)]
(\omega)
{{{}}}
\le \alpha 
\}.
\label{eq19} \end{align}
Thus, as our answer to the problem (G), we can assert the following theorem,
which is a slight generalization of our result in refs.
{{{}}}\cite{Ishi4}, \cite{Ishi8}.
\par
\noindent
\bf
Theorem 2
\rm
[Likelihood ratio test].
\sl
Assume the above notations.
Then,
the
${{{R} }}_{{{H_N}}}^{\epsilon{(\alpha)}}$
satisfies the condition (F).
And thus,
the rejection region
${\widehat R}^\alpha_{{{H_N}}}$
is given by
${{{R}}}_{{{H_N}}}^{{\epsilon(\alpha)} }$.
\par
\vskip0.5cm
\rm
\par
We believe that this theorem is the most orthodox answer to Problem (G).
However, in Section 2.2, we will propose another answer to Problem (G).

\vskip1.0cm
\par
\par
\section{The reverse relation between confidence interval method
and
statistical hypothesis testing
}
\rm
In this main section, we propose a new formulation of the confidence interval methods
and statistical hypothesis testing,
and show that they can be understood as two sides of the same coin
\subsection{Confidence interval method}
\par
\noindent
\par
\par
\noindent
\par
Let
${\mathsf O} = ({}X, {\cal F} , F{}){}$
be an observable
formulated in a
commutative $C^*$-algebra
${C_0(\Omega)}$.
Let $\Theta$ be a locally compact space with the 
semi-distance $d^x_{\Theta}$
$(\forall x \in X)$,
that is,
for each $x\in X$,
the map
$d^x_{\Theta}: \Theta^2 \to [0,\infty)$
satisfies that
(i):$d^x_\Theta (\theta, \theta )=0$,
(ii):$d^x_\Theta (\theta_1, \theta_2 )$
$=d^x_\Theta (\theta_2, \theta_1 )$,
(ii):$d^x_\Theta (\theta_1, \theta_3 )$
$\le d^x_\Theta (\theta_1, \theta_2 )
+
d^x_\Theta (\theta_2, \theta_3 )
$.
\noindent
\par
Let
$\pi:\Omega \to \Theta$
be a continuous map,
which is a kind of causal relation (in Axiom 2), and called
\it
{\lq\lq}quantity{\rq\rq}$.\; \;$
\rm
Let
$E:X \to \Theta$
be a continuous (or more generally, measurable) map,
which is called
\it
{\lq\lq}estimator{\rq\rq}$.\; \;$
\rm
\par
\noindent
\bf
Theorem 3
\rm
[Confidence interval method({\it cf.} ref.
{{{}}}{\cite{Ishi9})].
\sl
Let
$\gamma$
be a real number such that
$0 \ll \gamma < 1$,
for example,
$\gamma = 0.95$.
For any state
$ \omega ({}\in \Omega)$,
define
the positive number
$\eta^\gamma_{\omega}$
$({}> 0)$
such that:
\begin{align}
\eta^\gamma_{\omega}
=
\inf
\{
\eta > 0:
[F(\{ x \in X \;:\; 
d^x_\Theta ( E(x) , \pi( \omega ) )
< \eta
\}
)](\omega )
\ge \gamma
\}
\label{eq20}
\end{align}
Then we say that:
\rm
\begin{enumerate}
\item[(H$_1$)]
\it
\sl
the probability,
that
the measured value $x$
obtained
by the measurement
${\mathsf M}_{C_0(\Omega)} \big({}{\mathsf O}:= ({}X, {\cal F} , F{})  ,$
$ S_{[\omega_0 {}] } \big)$
satisfies the following
condition (\ref{eq21}),
is more than or equal to
$\gamma$
({}e.g., $\gamma= 0.95${}).
\end{enumerate}
\begin{align}
d^x_\Theta (E(x),  \pi(\omega_0){}) < {\eta }^\gamma_{\omega_0}
\label{eq21}
\end{align}
\sl
\par
\noindent
And further,
put
\begin{align}
D_{x}^{\gamma}
=
\{
\pi(\omega)
(\in
\Theta)
:
d^x_\Theta ({}E(x),
\pi(\omega )
)
<
\eta^\gamma_{\omega }
\}.
\label{eq22} 
\end{align}
which is called
\it
the $({}\gamma{})$-confidence interval.
\sl
Here,
we see the following equivalence:
\begin{align}
(\ref{eq21}) \; \Longleftrightarrow \;
\;
D_{x}^\gamma
\ni
\pi (\omega_0).
\label{eq23} 
\end{align}

\rm
\par
\noindent
\begin{center}
\unitlength=0.4mm
\begin{picture}(230,75)
\put(-19,0){{
\put(40,16){\scriptsize $x_0$}
\qbezier(40,20)(100,61)(157,42)
\qbezier(153,32)(200,-5)(257,32)
\path(107,49)(115,48)(107,45)
\put(112,53){$E$}
\put(105,-35){
\put(102,55){$\pi$}
\path(107,51)(99,48)(107,45)}
\put(40,20){\circle*{1}}
\put(157,41){\circle*{1}}
\put(155,45){\scriptsize $E({}x_0)$}
\put(151,33){\scriptsize $\; \pi(\omega_0)$}
\put(251,30){\scriptsize $\;\;\; \cdot \; \omega_0$}
\put(149,34){ \circle*{1} }
\put(175,35){$ D_{x_0}^\gamma$}
\put(153,63){ $\Theta$}
\put(253,63){ $\Omega$}
\put(57,63){$X$}
\allinethickness{0.5mm}
\put(60,30){\oval(70,60)}
\put(160,30){\oval(70,60)}
\put(260,30){\oval(70,60)}
\allinethickness{0.3mm}
\put(157,42){\ellipse{30}{30}}
}}
\put(40,-20){\bf Figure 3.
\rm
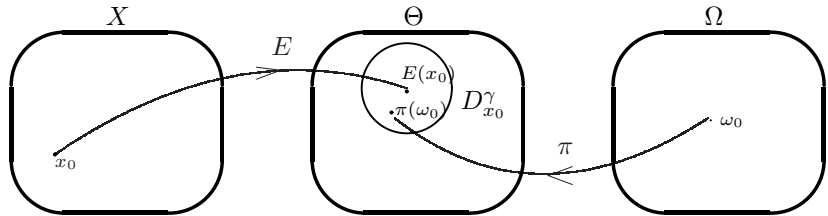
Confidence interval
$D^\gamma_{x_0}$
}
\end{picture}
\end{center}
\vskip1.5cm

\par
The following corollary 1 may not be useful.
However, it should be compared with Theorem 4.
\par
\noindent
\bf
Corollary 1
\sl
Further, consider a subset
$H_S
$
of $\Theta$,
which is called a 
"sure hypothesis".
Put
\begin{align}
{\widehat D}_{H_S}^{\gamma}
=
\bigcup_{\omega \in  \Omega \mbox{ \footnotesize such that }
\pi(\omega) \in {H_S}}
\{
E({x})
(\in
\Theta)
:
d^x_\Theta ({}E(x),
\pi(\omega )
)
<
\eta^\gamma_{\omega }
\}.
\label{eq126} 
\end{align}
\rm
\sl
Then we say that:
\rm
\begin{enumerate}
\item[(H$_2$)]
\it
\sl
the probability,
that
the measured value $x$
obtained
by the measurement
${\mathsf M}_{C_0(\Omega)} \big({}{\mathsf O}:= ({}X, {\cal F} , F{})  ,$
$ S_{[\ast {}] } (\pi^{-1}(H_S ) )\big)$
({\it cf.} (\ref{eq100}))
satisfies the following
condition (\ref{eq127}),
is more than or equal to
$\gamma$
({}e.g., $\gamma= 0.95${}).
\end{enumerate}
\begin{align}
{\widehat D}_{H_S}^\gamma
\ni
E(x).
\label{eq127} 
\end{align}

\par
\vskip1.0cm
\par
\rm
\subsection{Statistical hypothesis testing}
\par
\noindent
\par
The following theorem is our main theorem in this paper,
which says that it is contrary to Theorem 3
(the confidence interval method). In other words,they are two sides of the same coin.
\rm
\par
\noindent
\bf
Theorem 4
\rm
[Statistical hypothesis testing].
\sl
Let
$\alpha$
be a real number such that
$0 < \alpha \ll 1$,
for example,
$\alpha = 0.05$.
For any state
$ \omega ({}\in \Omega)$,
define
the positive number
$\eta^\alpha_{\omega}$
$({}> 0)$
such that:
\begin{align}
\eta^\alpha_{\omega}
&
=
\inf
\{
\eta > 0:
[F(\{ x \in X \;:\; 
d^x_\Theta ( E(x) , \pi( \omega ) )
\ge \eta
\}
)](\omega )
\le \alpha
\}
\nonumber
\\
\Big(
&=
\inf
\{
\eta > 0:
[F(\{ x \in X \;:\; 
d^x_\Theta ( E(x) , \pi( \omega ) )
< \eta
\}
)](\omega )
\ge 1- \alpha
\}
=\mbox{"$\eta^{1-\alpha}_\omega$ in (\ref{eq20})}"
\Big)
\label{eq24}
\end{align}
Then we say that:
\rm
\begin{enumerate}
\item[(I$_1$)]
\it
\sl
the probability,
that
the measured value $x$
obtained
by the measurement
${\mathsf M}_{C_0(\Omega)} \big({}{\mathsf O}:= ({}X, {\cal F} , F{})  ,$
$ S_{[\omega_0 {}] } \big)$
satisfies the following
condition (\ref{eq25}),
is less than or equal to
$\alpha$
({}e.g., $\alpha= 0.05${}).
\begin{align}
d^x_\Theta (E(x),  \pi(\omega_0){}) \ge  {\eta }^\alpha_{\omega_0}  .
\label{eq25}
\end{align}
\end{enumerate}
\vskip0.9cm

\sl
\par
\noindent
Further, consider a subset
$H_N
$
of $\Theta$,
which is called a 
"null hypothesis".
Put
\begin{align}
{\widehat R}_{H_N}^{\alpha}
=
\bigcap_{\omega \in  \Omega \mbox{ \footnotesize such that }
\pi(\omega) \in {H_N}}
\{
E({x})
(\in
\Theta)
:
d^x_\Theta ({}E(x),
\pi(\omega )
)
\ge
\eta^\alpha_{\omega }
\}.
\label{eq26} 
\end{align}
which is called
\it
the $({}\alpha{})$-rejection region
of
the null hypothesis
${H_N}$.
\rm
Then we say that:
\rm
\begin{enumerate}
\item[(I$_2$)]
\it
\sl
the probability,
that
the measured value $x$
obtained
by the measurement
${\mathsf M}_{C_0(\Omega)} \big({}{\mathsf O}:= ({}X, {\cal F} , F{})  ,$
$ S_{[\ast {}] } (\pi^{-1}( H_N )\big)$
({\it cf.} (\ref{eq100}))
satisfies the following
condition (\ref{eq27}),
is less than or equal to
$\alpha$
({}e.g., $\alpha= 0.05${}).
\end{enumerate}
\begin{align}
{\widehat R}_{H_N}^\alpha
\ni
E(x).
\label{eq27} 
\end{align}

\par
\noindent
\begin{center}
\unitlength=0.4mm
\begin{picture}(230,75)
\put(-19,0){{
\put(40,16){\scriptsize $x_0$}
\qbezier(40,20)(100,61)(157,42)
\qbezier(155,22)(200,-5)(257,32)
\path(107,49)(115,48)(107,45)
\put(112,53){$E$}
\put(105,-35){
\put(102,55){$\pi$}
\path(107,51)(99,46)(107,45)}
\put(40,20){\circle*{1}}
\put(157,41){\circle*{1}}
\put(155,45){\scriptsize $E({}x_0)$}
\put(151,25){\scriptsize $\; \pi(\omega_0)$}
\put(251,30){\scriptsize $\;\;\; \cdot \; \omega_0$}
\put(150,23){ \circle*{1} }
\put(175,35){$ {\widehat R}_{H_N}^\alpha$}
\put(153,63){ $\Theta$}
\put(253,63){ $\Omega$}
\put(57,63){$X$}
\allinethickness{0.5mm}
\put(60,30){\oval(70,60)}
\put(160,30){\oval(70,60)}
\put(260,30){\oval(70,60)}
\allinethickness{0.3mm}
\put(157,22){\ellipse{30}{30}}
}}
\put(40,-20){\bf Figure 4.
\rm
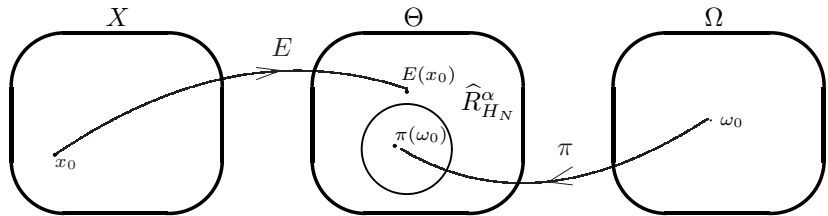
Rejection region
${\widehat R}^\alpha_{H_N}$
(when $H_N=\{\pi(\omega_0)\}$
}
\end{picture}
\end{center}

\par
\vskip1.0cm
\par

\par
\noindent
\bf
Remark 1
\rm
[The statistical meaning of Theorems 3 and 4].
(i):
The ${\widehat D}_{H_S}^\gamma$ in (\ref{eq126}) is the compliment of
${\widehat R}_{H_S}^\gamma$, however, Corollary 1 may not be useful.
\par
\noindent
(ii):
Consider the simultaneous measurement
${\mathsf M}_{C_0(\Omega)} \big({}{\mathsf O}^J:= ({}X^J, {\cal F}^J , F^J{})  ,$
$ S_{[\omega_0 {}] } \big)$,
and assume that a measured value $x=(x_1,x_2, \ldots , x_J)( \in X^J)$ is obtained by the simultaneous measurement.
Recall the formula (\ref{eq23}).
Then, it surely holds that
\begin{align}
\lim_{J \to \infty }
\frac{\mbox{Num} [\{ j \;|\; D_{x_j}^{\gamma} \ni \pi( \omega_0)]}{J}
\ge \gamma (= 0.95)
\label{eq28}
\end{align}
where
$\mbox{Num} [A]$ is the number of the elements of the set $A$.
Hence Theorem 3 can be tested by numerical analysis
(with random number). Similarly, Theorem 4 can be tested.

\par
\section{Examples}
The arguments in this section are continued from Example 2.
Let
$\alpha$
be a real number such that
$0 < \alpha \ll 1$,
for example,
$\alpha = 0.05$.
From the reverse relation between Theorem 3 (the confidence interval meyhod) 
and Theorem 4 (ststistical hypothesis testing),
Examples 4-10 in this section may be essentially the same as the examples of ref.\cite{Ishi9}.

\subsection{Population mean}
\par
\noindent
\rm
\par

\par
\noindent
\bf
Example 4
\rm
[Rejection region of $H_N=\{\mu_0\} \subseteq \Theta = {\mathbb R}$].
Consider the simultaneous measurement
${\mathsf M}_{C_0({\mathbb R} \times {\mathbb R}_+)}$
$({\mathsf O}_N^n = ({\mathbb R}^n, {\mathcal B}_{\mathbb R}^n, {{{N}}^n}) ,$
$S_{[(\mu, \sigma)]})$
in $C_0({\mathbb R} \times {\mathbb R}_+)$.
Thus,
we consider that
$\Omega = {\mathbb R} \times {\mathbb R}_+$,
$X={\mathbb R}^n$.
Assume that the real $\sigma$ in a state $\omega = (\mu, \sigma ) \in \Omega $
is fixed and known.
Put
$$
\Theta = {\mathbb R}
$$
\rm
The formula (\ref{eq16}) urges us to
define the estimator 
$E: {\mathbb R}^n \to \Theta  (\equiv {\mathbb R} )$
such that
\begin{align}
E(x)=E(x_1, x_2, \ldots , x_n )
=
\overline{\mu}(x)
=
\frac{x_1 + x_2 + \cdots + x_n}{n}
\label{eq29}
\end{align}
And consider the quantity $\pi: \Omega \to \Theta $ such that
$$
\Omega={\mathbb R} \times {\mathbb R}_+
\ni
\omega = (\mu, \sigma )
\mapsto \pi (\omega ) = \mu \in \Theta={\mathbb R}
$$
Consider the following semi-distance $d_{\Theta}^{(1)}$
in $\Theta (={\mathbb R} )$:
\begin{align}
d_{\Theta}^{(1)}(\theta_1, \theta_2)
=
|\theta_1 - \theta_2|
\label{eq30}
\end{align}
Define the null hypothesis $H_N$
such that
$$
H_N=\{\mu_0\} (\subseteq \Theta (= {\mathbb R}))
$$
For any
$ \omega=(\mu, \sigma )  ({}\in \Omega=
{\mathbb  R} \times {\mathbb R}_+ )$,
define
the positive number
$\eta^\alpha_{\omega}$
$({}> 0)$
such that:
\begin{align}
\eta^\alpha_{\omega}
=
\sup
\{
\eta > 0:
[F ({}E^{-1} ({}
{ {\rm Ball}^C_{d_{\Theta}^{(1)}}}(\pi(\omega) ; \eta{}))](\omega )
\le \alpha
\}
\nonumber
\end{align}
where
${{\rm Ball}^C_{d_{\Theta}^{(1)}}}(\pi( \omega ) ; \eta)$
$=$
$\{ \theta
({}\in\Theta):
d_{\Theta}^{(1)} ({}\mu, \theta {}) \ge \eta \}$
$=
\Big(
( -\infty, \mu - \eta] \cup  [\mu + \eta , \infty )
\Big)
$

Hence we see that
\begin{align}
&
E^{-1}({{\rm Ball}^C_{d_{\Theta}^{(1)}}}(\pi (\omega) ; \eta ))
=
E^{-1}
\Big(
( -\infty, \mu - \eta] \cup  [\mu + \eta , \infty )
\Big)
\nonumber
\\
=
&
\{
(x_1, \ldots , x_n )
\in {\mathbb R}^n
\;:
\;
\frac{x_1+\ldots + x_n }{n} 
\le
\mu - \eta 
\mbox{ or }
\mu + \eta 
\le
\frac{x_1+\ldots + x_n }{n} 
\}
\nonumber
\\
=
&
\{
(x_1, \ldots , x_n )
\in {\mathbb R}^n
\;:
\;
|\frac{(x_1- \mu)+\ldots + (x_n- \mu) }{n} 
|
\ge \eta
\}
\label{eq31}
\end{align}
Thus,
\begin{align}
&
[{{{N}}}^n
(E^{-1}({{\rm Ball}^C_{d_{\Theta}^{(1)}}}(\pi(\omega) ; \eta ))]
({}\omega{})
\nonumber
\\
=
&
\frac{1}{({{\sqrt{2 \pi }\sigma{}}})^n}
\underset{{
|\frac{(x_1- \mu)+\ldots + (x_n- \mu) }{n} 
|
\ge \eta
}}{\int \cdots \int}
\exp[{}- \frac{\sum_{k=1}^n ({}{}{x_k} - {}{\mu}  {})^2 
}
{2 \sigma^2}    {}] d {}{x_1} d {}{x_2}\cdots dx_n
\nonumber
\\
=
&
\frac{1}{({{\sqrt{2 \pi }\sigma{}}})^n}
\underset{{
|
\frac{x_1+\ldots + x_n }{n} 
|
\ge   \eta 
}}{\int \cdots \int}
\exp[{}- \frac{\sum_{k=1}^n ({}{}{x_k}  {}{}  {})^2 
}
{2 \sigma^2}    {}] d {}{x_1} d {}{x_2}\cdots dx_n
\nonumber
\\
=
&
\frac{\sqrt{n}}{{\sqrt{2 \pi }\sigma{}}}
\int_{{x \ge \eta}} \exp[{}- \frac{{n}{x}^2 }{2 \sigma^2}] d {x}
=
\frac{1}{{\sqrt{2 \pi }{}}}
\int_{{x \ge  \sqrt{n} \eta/\sigma}} \exp[{}- \frac{{x}^2 }{2 }] d {x}
\label{eq32}
\end{align}
Solving the following equation:
\begin{align}
\frac{1}{{\sqrt{2 \pi }{}}}
\int^{-z(\alpha/2)}_{-\infty} \exp[{}- \frac{{x}^2 }{2 }] d {x}
=
\frac{1}{{\sqrt{2 \pi }{}}}
\int_{z(\alpha/2)}^{\infty} \exp[{}- \frac{{x}^2 }{2 }] d {x}
=
\frac{\alpha}{2}
\label{eq33}
\end{align}
we define that
\begin{align}
\eta^\alpha_{\omega} = 
 \frac{\sigma}{\sqrt{n}}
 z(\frac{\alpha}{2})
\label{eq34}
\end{align}

Therefore,
we get ${\widehat R}_{H_N}^{\alpha}$
(
the $({}\alpha{})$-rejection region
of
$H_N
(= \{ \mu_0\} \subseteq \Theta (= {\mathbb R}))$
)
as follows:
\begin{align}
{\widehat R}_{\{ \mu_0 \}}^{\alpha}
&
=
\bigcap_{\pi(\omega ) =\mu \in \{ \mu_0 \}  }
\{
{E(x)}
(\in
\Theta=
{\mathbb R})
:
d_{\Theta}^{(1)} ({}E(x),
\pi (\omega))
\ge
\eta^\alpha_{\omega }
\}
\nonumber
\\
&
=
\{ E(x) (= \frac{x_1+ \ldots + x_n}{n}) \in {\mathbb R} 
\;:\;
\overline{\mu}(x)
- \mu_0
=
\frac{x_1+ \ldots + x_n}{n}
-
\mu_0
\ge 
 \frac{\sigma}{\sqrt{n}}
 z(\frac{\alpha}{2})
 \}
\label{eq35}
\end{align}
\bf
Remark 2
\rm
Note that
the ${\widehat R}_{\{\mu_0\}}^{\alpha}$
(
the $({}\alpha{})$-rejection region
of
$\{ \mu_0\}$
)
depends on $\sigma$.
%
\rm
Thus, putting
\begin{align}
{\widehat R}_{\{ \mu_0 \} \times {\mathbb R}_+}^{\alpha}
=
\{  (\overline{\mu}(x), \sigma) \in {\mathbb R} \times {\mathbb R}_+
\;:\;
| \mu_0 - \overline{\mu}(x)|
=
| \mu_0 - \frac{x_1+ \ldots + x_n}{n}|
\ge 
 \frac{\sigma}{\sqrt{n}}
 z(\frac{\alpha}{2})
 \}
\label{eq36}
\end{align}
we see that
${\widehat R}_{ \{\mu_0\}\times {\mathbb R}_+}^{\alpha}$="the slash part in Figure 5".

\par
\noindent
\begin{center}
\unitlength=0.43mm
\begin{picture}(250,150)
\put(75,50)
{{
\put(100,-10){${\mathbb R}$}
\put(-8,85){$\sigma$}
\put(80,20){\large{${\widehat R}_{\{\mu_0\} \times {\mathbb R}_+}^{\alpha}$}}
\put(-30,0){\vector(1,0){130}}
\put(-15,0){\vector(0,1){90}}
\allinethickness{0.1mm}
\multiput(-43,58)(3,-3){20}{\line(-1,0){30}}
\multiput(74,58)(-3,-3){20}{\line(1,0){30}}
\path(15,-2)(15,2)
\put(14,-8){${\mu}_0$}
\put(25,50){}
\thicklines
\put(-50,65){\line(1,-1){65}}
\put(81,65){\line(-1,-1){65}}
}}
\put(10,20){\bf Figure 5.
\rm
Rejection region
${\widehat R}_{\{\mu_0\}}^{\alpha}$
(which depends on $\sigma$)
%
}
\end{picture}
\end{center}
%

\par
\noindent
\bf
Example 5
\rm
[Rejection region of $H_N=( -\infty , \mu_0] \subseteq \Theta (={\mathbb R})$]. Consider the simultaneous measurement
${\mathsf M}_{C_0({\mathbb R} \times {\mathbb R}_+)}$
$({\mathsf O}_N^n = ({\mathbb R}^n, {\mathcal B}_{\mathbb R}^n, {{{N}}^n}) ,$
$S_{[(\mu, \sigma)]})$
in $C_0({\mathbb R} \times {\mathbb R}_+)$.
Thus,
we consider that
$\Omega = {\mathbb R} \times {\mathbb R}$,
$X={\mathbb R}^n$.
Assume that the real $\sigma$ in a state $\omega = (\mu, \sigma ) \in \Omega $
is fixed and known.
Put
$$
\Theta = {\mathbb R}
$$
\rm
The formula (\ref{eq16}) urges us to
define the estimator 
$E: {\mathbb R}^n \to \Theta  (\equiv {\mathbb R} )$
such that
\begin{align}
E(x)=
=
\overline{\mu}(x)
=
\frac{x_1 + x_2 + \cdots + x_n}{n}
\label{eq37}
\end{align}
And consider the quantity $\pi: \Omega \to \Theta $ such that
$$
\Omega={\mathbb R} \times {\mathbb R}_+
\ni
\omega = (\mu, \sigma )
\mapsto \pi (\omega ) = \mu \in \Theta={\mathbb R}
$$
Consider the following semi-distance $d_{\Theta}^{(2)}$
in $\Theta (={\mathbb R} )$:
\begin{align}
d_{\Theta}^{(2)}((\theta_1, \theta_2)
=
\cases
|\theta_1 - \theta_2| \quad & \theta_0 \le \theta_1, \theta_2 
\\
|\theta_2 - \theta_0| \quad & \theta_1 \le \theta_0 \le \theta_2 
\\
|\theta_1 - \theta_0| \quad & \theta_2 \le \theta_0 \le \theta_1
\\
0 \quad & \theta_1  , \theta_2 \le \theta_0
\endcases
\label{eq38}
\end{align}
Define the null hypothesis $H_N$
such that
$$
H_N= ( -\infty , \mu_0] (\subseteq \Theta (= {\mathbb R}))
$$
For any
$ \omega=(\mu, \sigma )  ({}\in \Omega=
{\mathbb  R} \times {\mathbb R}_+ )$,
define
the positive number
$\eta^\alpha_{\omega}$
$({}> 0)$
such that:
\begin{align}
\eta^\alpha_{\omega}
=
\sup
\{
\eta > 0:
[F ({}E^{-1} ({}
{ {\rm Ball}^C_{d_{\Theta}^{(2)}}}(\pi(\omega) ; \eta{}))](\omega )
\le \alpha
\}
\nonumber
\end{align}
where
${{\rm Ball}^C_{d_{\Theta}^{(2)}}}(\pi( \omega ) ; \eta)$
$=$
$\{ \theta
({}\in\Theta):
d_{\Theta}^{(2)} ({}\mu, \theta {}) \ge \eta \}$
$=
\Big(
( -\infty, \mu - \eta] \cup  [\mu + \eta , \infty )
\Big)
$

Hence we see that
\begin{align}
&
E^{-1}({{\rm Ball}^C_{d_{\Theta}^{(2)}}}(\pi (\omega) ; \eta ))
=
E^{-1}
\Big(
  [\mu + \eta , \infty )
\Big)
\nonumber
\\
=
&
\{
(x_1, \ldots , x_n )
\in {\mathbb R}^n
\;:
\;
\mu + \eta 
\le
\frac{x_1+\ldots + x_n }{n} 
\}
\nonumber
\\
=
&
\{
(x_1, \ldots , x_n )
\in {\mathbb R}^n
\;:
\;
\frac{(x_1- \mu)+\ldots + (x_n- \mu) }{n} 
\ge \eta
\}
\label{eq39}
\end{align}
Thus,
\begin{align}
&
[{{{N}}}^n
(E^{-1}({{\rm Ball}^C_{d_{\Theta}^{(2)}}}(\pi(\omega) ; \eta ))]
({}\omega{})
\nonumber
\\
=
&
\frac{1}{({{\sqrt{2 \pi }\sigma{}}})^n}
\underset{{
\frac{(x_1- \mu)+\ldots + (x_n- \mu) }{n} 
\ge \eta
}}{\int \cdots \int}
\exp[{}- \frac{\sum_{k=1}^n ({}{}{x_k} - {}{\mu}  {})^2 
}
{2 \sigma^2}    {}] d {}{x_1} d {}{x_2}\cdots dx_n
\nonumber
\\
=
&
\frac{1}{({{\sqrt{2 \pi }\sigma{}}})^n}
\underset{{
\frac{x_1+\ldots + x_n }{n} 
\ge   \eta 
}}{\int \cdots \int}
\exp[{}- \frac{\sum_{k=1}^n ({}{}{x_k}  {}{}  {})^2 
}
{2 \sigma^2}    {}] d {}{x_1} d {}{x_2}\cdots dx_n
\nonumber
\\
=
&
\frac{\sqrt{n}}{{\sqrt{2 \pi }\sigma{}}}
\int_{{|x| \ge \eta}} \exp[{}- \frac{{n}{x}^2 }{2 \sigma^2}] d {x}
=
\frac{1}{{\sqrt{2 \pi }{}}}
\int_{{|x| \ge  \sqrt{n} \eta/\sigma}} \exp[{}- \frac{{x}^2 }{2 }] d {x}
\label{eq40}
\end{align}
Solving the following equation:
\begin{align}
\frac{1}{{\sqrt{2 \pi }{}}}
\int^{-z(\alpha/2)}_{-\infty} \exp[{}- \frac{{x}^2 }{2 }] d {x}
=
\frac{1}{{\sqrt{2 \pi }{}}}
\int_{z(\alpha/2)}^{\infty} \exp[{}- \frac{{x}^2 }{2 }] d {x}
=
{\alpha}
\label{eq41}
\end{align}
we define that
\begin{align}
\eta^\alpha_{\omega} = 
 \frac{\sigma}{\sqrt{n}}
 z({\alpha})
\label{eq42}
\end{align}

Therefore,
we get ${\widehat R}_{H_N}^{\alpha}$
(
the $({}\alpha{})$-rejection region
of
$H_N
(= ( - \infty , \mu_0] \subseteq \Theta (= {\mathbb R}))$
)
as follows:
\begin{align}
{\widehat R}_{ ( - \infty , \mu_0]}^{\alpha}
&
=
\bigcap_{\pi(\omega ) = \mu \in  ( - \infty , \mu_0] }
\{
{E(x)}
(\in
\Theta=
{\mathbb R})
:
d_{\Theta}^{(2)} ({}E(x),
\pi (\omega))
\ge
\eta^\alpha_{\omega }
\}
\nonumber
\\
&
=
\{ E(x) (= \frac{x_1+ \ldots + x_n}{n}) \in {\mathbb R} 
\;:\;
 \frac{x_1+ \ldots + x_n}{n} - \mu_0
\ge 
 \frac{\sigma}{\sqrt{n}}
 z({\alpha})
 \}
\label{eq43}
\end{align}
\rm
Thus, in a similar way of Remark 2, we see that ${\widehat R}_{
(- \infty, \mu_0 ] \times {\mathbb R}_+}^{\alpha}$="the slash part in Figure 6", where
\begin{align}
{\widehat R}_{( - \infty, \mu_0 ] \times {\mathbb R}_+}^{\alpha}
=
\{ (E(x) (= \frac{x_1+ \ldots + x_n}{n}), \sigma) \in {\mathbb R} \times {\mathbb R}_+
\;:\;
\frac{x_1+ \ldots + x_n}{n}
 -
 \mu_0
\ge 
 \frac{\sigma}{\sqrt{n}}
 z({\alpha})
 \}
\label{eq44}
\end{align}

\par
\noindent
\begin{center}
\unitlength=0.43mm
\begin{picture}(250,150)
\put(75,50)
{{
\put(100,-10){${\mathbb R}$}
\put(-10,85){$\sigma$}
\put(100,10){\large{${\widehat R}_{(-\infty,\mu_0] \times {\mathbb R}_+}^{\alpha}$}}
\put(-30,0){\vector(1,0){130}}
\put(-15,0){\vector(0,1){90}}
\allinethickness{0.1mm}
\path(15,-2)(15,2)
\put(14,-8){${\mu}_0$}
\put(25,50){}
\put(45,0){
\multiput(74,58)(-5,-3){20}{\line(1,0){30}}
\thicklines
\put(82,65){\line(-5,-3){110}}
}
}}
\put(10,20){\bf Figure 6.
\rm
Rejection region
${\widehat R}_{(-\infty,\mu_0]}^{\alpha}$
(which depends on $\sigma$)
%
}
\end{picture}
\end{center}
%


\par
\noindent
\par
\subsection{Population variance}
\par
\noindent
\bf
Example 6
\rm
[Rejection region of $H_N= \{ \sigma_0 \} \subseteq \Theta (={\mathbb R}_+ $].
Consider the simultaneous measurement
${\mathsf M}_{C_0({\mathbb R} \times {\mathbb R}_+)}$
$({\mathsf O}_N^n = ({\mathbb R}^n, {\mathcal B}_{\mathbb R}^n, {{{N}}^n}) ,$
$S_{[(\mu, \sigma)]})$
in $C_0({\mathbb R} \times {\mathbb R}_+)$.
Thus,
we consider that
$\Omega = {\mathbb R} \times {\mathbb R}_+$,
$X={\mathbb R}^n$.
Assume that the real $\mu$ in a state $\omega = (\mu, \sigma ) \in \Omega $
is fixed and known.
Put
$$
\Theta = {\mathbb R}_+
$$
\rm
The formula (\ref{eq16}) may urge us to
define the estimator 
$E: {\mathbb R}^n \to \Theta  (\equiv {\mathbb R}_+ )$
such that
\begin{align}
E(x)=E(x_1, x_2, \ldots , x_n )
=
\overline{\sigma}(x)=
\sqrt{
\frac{\sum_{k=1}^n ( x_k - \overline{\mu}(x))^2}{n}
}
\label{eq45}
\end{align}
And consider the quantity $\pi: \Omega \to \Theta $ such that
$$
\Omega={\mathbb R} \times {\mathbb R}_+
\ni
\omega = (\mu, \sigma )
\mapsto \pi (\omega ) = \sigma \in \Theta={\mathbb R}_+
$$
Define the null hypothesis $H_N$
such that
$$
H_N=\{\sigma_0\} (\subseteq \Theta (= {\mathbb R}_+))
$$
Consider the following semi-distance $d_{\Theta}^{(1)}$
in $\Theta (={\mathbb R}_+ )$:
\begin{align}
d_{\Theta}^{(1)}(\theta_1, \theta_2)
=
|
\int_{\sigma_1}^{\sigma_2} \frac{1}{\sigma} d \sigma
|
=
|\log{\sigma_1} - \log{\sigma_2} |
\label{eq46}
\end{align}

%
%
%

For any
$ \omega=(\mu, {\sigma} )  ({}\in
\Omega=
{\mathbb  R} \times {\mathbb R}_+ )$,
define
the positive number
$\eta^\alpha_{\omega}$
$({}> 0)$
such that:
\begin{align}
\eta^\alpha_{\omega}
=
\sup
\{
\eta > 0:
[F ({}E^{-1} ({}
{{\rm Ball}^C_{d_{\Theta}^{(1)}}}(\omega ; \eta{}))](\omega )
\le \alpha
\}
\label{eq47}
\end{align}
where
\begin{align}
{{\rm Ball}^C_{d_{\Theta}^{(1)}}}(\omega ; \eta )
=
{{\rm Ball}^C_{d_{\Theta}^{(1)}}}((\mu ; {\sigma} ), \eta )
=
{\mathbb R} \times \{ \sigma' \;:\; |\log(\sigma'/{\sigma})| \ge \eta
\}
=
{\mathbb R} \times \big(
(0,{\sigma} e^{-\eta}] \cup  [{\sigma} e^{\eta} , \infty )
\big)
\label{eq48}
\end{align}
Then,
\begin{align}
&
E^{-1}( {{\rm Ball}^C_{d_{\Theta}^{(1)}}}(\omega ; \eta ))
=
E^{-1}
\Big(
{\mathbb R} \times \big(
(0,{\sigma} e^{-\eta}] \cup  [{\sigma} e^{\eta} , \infty )
\big)
\Big)
\nonumber
\\
=
&
\{
(x_1, \ldots , x_n )
\in
{\mathbb R}^n
\;:
\;
\Big(
\frac{\sum_{k=1}^n ( x_k - 
\overline{\mu}
(x))^2}{n}
\Big)^{1/2}
\le
{\sigma} e^{ -\eta }
\mbox{ or }
{\sigma} e^{ \eta }
\le
\Big(
\frac{\sum_{k=1}^n ( x_k - 
\overline{\mu}
(x))^2}{n}
\Big)^{1/2}
\}
\label{eq49}
\end{align}
Hence we see, by (\ref{eq10}), that
\begin{align}
&
[{{{N}}}^n
(E^{-1}({{\rm Ball}^C_{d_{\Theta}^{(1)}}}(\omega; \eta ))]
({}\omega{})
\nonumber
\\
=
&
\frac{1}{({{\sqrt{2 \pi }{\sigma}{}}})^n}
\underset{{
E^{-1}
\Big(
{\mathbb R} \times \big(
(0,{\sigma} e^{-\eta}] \cup  [{\sigma} e^{\eta} , \infty )
\big)
\Big)
}}{\int \cdots \int}
\exp[{}- \frac{\sum_{k=1}^n ({}{}{x_k} - {}{\mu}  {})^2 
}
{2 {\sigma}^2}    {}] d {}{x_1} d {}{x_2}\cdots dx_n
\nonumber
\\
=
&
\int_0^{{n}  e^{- 2 \eta}}
p^{\chi^2}_{n-1} (x ) 
dx
+
\int_{{n}  e^{ 2 \eta}}^\infty 
p^{\chi^2}_{n-1} (x ) 
dx
=
1-
\int_{{n}  e^{- 2 \eta}}^{{n}  e^{ 2 \eta}} 
p^{\chi^2}_{n-1} (x ) 
dx
\label{eq50}
\end{align}
Using the chi-squared distribution $p^{{\chi}^2}_{n-1}({ x} )$
(with $n-1$ degrees of freedom) in (\ref{eq11}),
define the $\eta^\alpha_{\omega}$ such that
\begin{align}
1-
\alpha
=
\int_{{n} e^{-2 \eta^\alpha_{\omega}}}^{{n}  e^{2 \eta^\alpha_{\omega}}} 
p^{\chi^2}_{n-1} (x ) dx
\label{eq51}
\end{align}
where it should be noted that
the $\eta^\alpha_{\omega}$
depends on only $\alpha$ and $n$.
Thus, put
\begin{align}
\eta^\alpha_{\omega} = \eta^\alpha_{n}
\label{eq52}
\end{align}
Hence we get
the
${\widehat R}_{H_N}^{\alpha}$
(
the $({}\alpha{})$-rejection region
of
$H_N= \{ \sigma_0 \} \subseteq \Theta ={\mathbb R}_+$
)
as follows:
\begin{align}
{\widehat R}_{H_N}^{\alpha}
&
=
{\widehat R}_{\{ \sigma_0 \}}^{\alpha}
=
\bigcap_{\pi(\omega ) = \sigma \in  \{\sigma_0 \}}
\{
{E(x)}
(\in
\Theta)
:
d^{(2)}_\Theta  ({}E(x),
\omega)
\ge
\eta^\alpha_{\omega}
\}
\nonumber
\\
&
=
\{
{E(x)}
(\in
\Theta={\mathbb R}_+)
:
d^{(2)}_\Theta  ({}E(x),
(\mu, \sigma_0))
\ge
\eta^\alpha_{n}
\}
\nonumber
\\
&
=
\{ \overline{\sigma}(x) (\in \Theta ={\mathbb R}_+ )
\;:
\;
\overline{\sigma}(x)
\le
{\sigma_0} e^{ -\eta^\alpha_{n} }
\mbox{ or }
{\sigma_0} e^{\eta^\alpha_{n} }
\le
\overline{\sigma}(x)
\}
\label{eq53}
\end{align}
where
$\overline{\sigma}(x)
=
\Big(
\frac{\sum_{k=1}^n ( x_k - 
\overline{\mu}
(x))^2}{n}
\Big)^{1/2}
$.

Thus, in a similar way of Remark 2, we see that ${\widehat R}_{{\mathbb R} \times \{\sigma_0 \}}^{\alpha}$="the slash part in Figure 7", where
\begin{align}
{\widehat R}_{{\mathbb R} \times \{\sigma_0 \}}^{\alpha}
=
\{ (\mu , \overline{\sigma}(x) ) \in {\mathbb R} \times {\mathbb R}_+
\;:\;
\overline{\sigma}(x)
\le
{\sigma_0} e^{ -\eta^\alpha_{n} }
\mbox{ or }
{\sigma_0} e^{\eta^\alpha_{n} }
\le
\overline{\sigma}(x)
\}
\label{eq54}
\end{align}

\par
\noindent
\begin{center}
\unitlength=0.43mm
\begin{picture}(250,150)
\put(75,50)
{{
\put(100,-9){$\mu$}
\put(-10,90){${\mathbb R}_+$}
\put(100,66){\large{${\widehat R}_{{\mathbb R} \times \{ \sigma_0 \}}^{\alpha}$}}
\put(-30,0){\vector(1,0){130}}
\put(-15,0){\vector(0,1){90}}
\put(25,50){}
\thicklines
\allinethickness{0.1mm}
\multiput(-30,57)(4,0){30}{\line(0,1){30}}
\multiput(-30,25)(4,0){30}{\line(0,-1){25}}
\put(-33,56){\line(1,0){135}}
\put(-33,25){\line(1,0){135}}
\put(-33,65){
\put(30,5){\vector(-1,-1){10}}
\put(30,5){
${{\sigma}_0
e^{\eta^\alpha_{n} }
}
$
}
}
\put(-33,35){
\put(30,5){\vector(-1,0){10}}
\put(32,3){$\sigma_0$}
}
\put(-33,5){
\put(30,5){\vector(-1,1){10}}
\put(30,5){
${{\sigma}_0
e^{- \eta^\alpha_{n}}
}
$
}
}
}}
\put(20,25){\bf Figure 7.
\rm
Rejection region
$
{\widehat R}_{\{ \sigma_0 \}}^{\alpha}
$
}
\end{picture}
\end{center}
%
%
%
%


\par
\noindent

\par
\noindent
\bf
Example 7
\rm
[Rejection region of $H_N= (-\infty, \sigma_0 ] \subseteq \Theta (={\mathbb R}_+ $].
Consider the simultaneous measurement
${\mathsf M}_{C_0({\mathbb R} \times {\mathbb R}_+)}$
$({\mathsf O}_N^n = ({\mathbb R}^n, {\mathcal B}_{\mathbb R}^n, {{{N}}^n}) ,$
$S_{[(\mu, \sigma)]})$
in $C_0({\mathbb R} \times {\mathbb R}_+)$.
Thus,
we consider that
$\Omega = {\mathbb R} \times {\mathbb R}_+$,
$X={\mathbb R}^n$.
Assume that the real $\mu$ in a state $\omega = (\mu, \sigma ) \in \Omega $
is fixed and known.
Put
$$
\Theta = {\mathbb R}_+
$$
\rm
The formula (\ref{eq16}) may urge us to
define the estimator 
$E: {\mathbb R}^n \to \Theta  (\equiv {\mathbb R}_+ )$
such that
\begin{align}
E(x)=E(x_1, x_2, \ldots , x_n )
=
\overline{\sigma}(x)=
\sqrt{
\frac{\sum_{k=1}^n ( x_k - \overline{\mu}(x))^2}{n}
}
\label{eq55}
\end{align}
And consider the quantity $\pi: \Omega \to \Theta $ such that
$$
\Omega={\mathbb R} \times {\mathbb R}_+
\ni
\omega = (\mu, \sigma )
\mapsto \pi (\omega ) = \sigma \in \Theta={\mathbb R}_+
$$
Define the null hypothesis $H_N$
such that
$$
H_N=( - \infty ,\sigma_0] (\subseteq \Theta (= {\mathbb R}_+))
$$
Consider the following semi-distance $d_{\Theta}^{(2)}$
in ${\mathbb R} \times {\mathbb R}_+$:
\begin{align}
d_{\Theta}^{(2)}((\mu_1,\sigma_1), (\mu_2,\sigma_2))
=
\cases
|
\int_{\sigma_1}^{\sigma_2} \frac{1}{\sigma} d \sigma
|
=
|\log{\sigma_1} - \log{\sigma_2}|
\quad & ( \sigma_0 \le \sigma_1, \sigma_2 )
\\
|
\int_{\sigma_0}^{\sigma_2} \frac{1}{\sigma} d \sigma
|
=
|\log{\sigma_0} - \log{\sigma_2}|
\quad &  ( \sigma_1 \le \sigma_0 \le \sigma_2 )
\\
|
\int_{\sigma_0}^{\sigma_1} \frac{1}{\sigma} d \sigma
|
=
|\log{\sigma_0} - \log{\sigma_1}|
\quad &  ( \sigma_2 \le \sigma_0 \le \sigma_1 )
\\
0
\quad & ( \sigma_1, \sigma_2 \le \sigma_0 )
\endcases
\label{eq56}
\end{align}
For any
$ \omega=(\mu, {\sigma} )  ({}\in
\Omega=
{\mathbb  R} \times {\mathbb R}_+ )$,
define
the positive number
$\eta^\alpha_{\omega}$
$({}> 0)$
such that:
\begin{align}
\eta^\alpha_{\omega}
=
\sup
\{
\eta > 0:
[F ({}E^{-1} ({}
{{\rm Ball}^C_{d_{\Theta}^{(2)}}}(\omega ; \eta{}))](\omega )
\le \alpha
\}
\label{eq57}
\end{align}
where
\begin{align}
{{\rm Ball}^C_{d_{\Theta}^{(2)}}}(\omega ; \eta )
=
{{\rm Ball}^C_{d_{\Theta}^{(2)}}}((\mu ; {\sigma} ), \eta )
=
{\mathbb R} \times [
 \sigma e^{\eta}, \infty 
)
\label{eq58}
\end{align}
Then,
\begin{align}
&
E^{-1}( {{\rm Ball}^C_{d_{\Theta}^{(2)}}}(\omega ; \eta ))
=
E^{-1}
\Big(
[
\sigma  e^{\eta}, \infty 
)
\Big)
\nonumber
\\
=
&
\{
(x_1, \ldots , x_n )
\in
{\mathbb R}^n
\;:
\;
\sigma e^{\eta}
\le
\overline{\sigma}(x)
=
\Big(
\frac{\sum_{k=1}^n ( x_k - 
\overline{\mu}
(x))^2}{n}
\Big)^{1/2}
\}
\label{eq59}
\end{align}
Hence we see, by (\ref{eq10}), that
\begin{align}
&
[{{{N}}}^n
(E^{-1}({{\rm Ball}^C_{d_{\Theta}^{(2)}}}(\omega; \eta ))]
({}\omega{})
\nonumber
\\
=
&
\frac{1}{({{\sqrt{2 \pi }{\sigma}{}}})^n}
\underset{{
\sigma_0 e^{\eta}
\le
\overline{\sigma}(x)
}}{\int \cdots \int}
\exp[{}- \frac{\sum_{k=1}^n ({}{}{x_k} - {}{\mu}  {})^2 
}
{2 {\sigma}^2}    {}] d {}{x_1} d {}{x_2}\cdots dx_n
\nonumber
\nonumber
\\
=
&
\int_{\frac{{n}  e^{ 2 \eta}\sigma^2}{\sigma^2}}^\infty 
p^{\chi^2}_{n-1} (x ) 
dx
\nonumber
\\
\le
&
\int_{{n}  e^{ 2 \eta}}^\infty 
p^{\chi^2}_{n-1} (x ) 
dx
\label{eq60}
\end{align}
Solving the following equation,
define the $(\eta^\alpha_{n})' (>0)$ such that
\begin{align}
\alpha
=
\int_{{n} e^{2 (\eta^\alpha_{n})'}}^{\infty} 
p^{\chi^2}_{n-1} (x ) dx
\label{eq61}
\end{align}
Hence we get
the
${\widehat R}_{H_N}^{\alpha}$
(
the $({}\alpha{})$-rejection region
of
$H_N= {\mathbb R} \times (0, \sigma_0] $
)
as follows:
\begin{align}
{\widehat R}_{H_N}^{\alpha}
&
=
{\widehat R}_{
{\mathbb R} \times (0, \sigma_0 ]
}^{\alpha}
=
\bigcap_{\pi(\omega ) = \omega \in {\mathbb R} \times (0, \sigma_0 ]}
\{
{E(x)}
(\in
\Omega)
:
d^{(2)}_\Theta  ({}E(x),
\omega)
\ge
\eta^\alpha_{\omega}
\}
\nonumber
\\
&
%
%
%
%
%
%
=
\{
{E(x)}
(\in
\Omega)
:
d^{(2)}_\Theta  ({}E(x),
\omega)
\ge
(\eta^\alpha_{n})'
\}
\nonumber
\\
&
=
\{ (\mu, {\sigma}(=
\overline{\sigma}(x)
) ) \in {\mathbb R} \times {\mathbb R}_+
\;:
\;
{\sigma_0} e^{(\eta^\alpha_{n})' }
\le
\overline{\sigma}(x)
\}
\label{eq62}
\end{align}
where
$\overline{\sigma}(x)
=
\Big(
\frac{\sum_{k=1}^n ( x_k - 
\overline{\mu}
(x))^2}{n}
\Big)^{1/2}
$.

Thus, in a similar way of Remark 2, we see that ${\widehat R}_{{\mathbb R} \times (0, \sigma_0]}^{\alpha}$="the slash part in Figure 8", where
\begin{align}
{\widehat R}_{{\mathbb R} \times (0, \sigma_0]}^{\alpha}
=
\{ (\mu , \overline{\sigma}(x) ) \in {\mathbb R} \times {\mathbb R}_+
\;:\;
\;:
\;
{\sigma_0} e^{(\eta^\alpha_{n})' }
\le
\overline{\sigma}(x)
\}
\label{eq63}
\end{align}

\par
\noindent
\begin{center}
\unitlength=0.43mm
\begin{picture}(250,150)
\put(75,50)
{{
\put(100,-8){$\mu$}
\put(-10,90){${\mathbb R}_+$}
\put(100,66){\large{${\widehat R}_{{\mathbb R} \times (0, \sigma_0]}^{\alpha}$}}
\put(-30,0){\vector(1,0){130}}
\put(-15,0){\vector(0,1){90}}
\put(25,50){}
\thicklines
\allinethickness{0.1mm}
\multiput(-30,57)(4,0){30}{\line(0,1){30}}
\put(-33,56){\line(1,0){135}}
\put(-33,25){\line(1,0){135}}
\put(-33,65){
\put(30,5){\vector(-1,-1){10}}
\put(30,5){
${{\sigma}_0
e^{(\eta^\alpha_{n})' }
}
$
}
}
\put(-33,35){
\put(30,5){\vector(-1,0){10}}
\put(32,3){$\sigma_0$}
}
\put(-33,5){
\put(30,5){\vector(-1,1){10}}
\put(30,5){
${{\sigma}_0
e^{- (\eta^\alpha_{n})'}
}
$
}
}
}}
\put(20,25){\bf Figure 8.
\rm
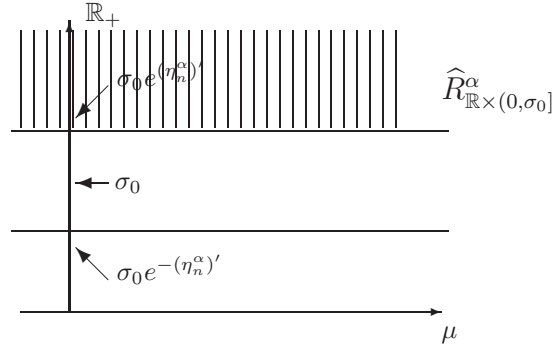
Rejection region
$
{\widehat R}_{ (0, \sigma_0]}^{\alpha}
$
}
\end{picture}
\end{center}

%
%
\subsection{The difference of the population means
}
\par
\noindent
\rm
\par
The arguments in this section are continued from Example 2.
\par
\noindent
\bf
Example 8
\rm
\rm
[Rejection region in the case that "$\pi(\mu_1, \mu_2)=\mu_1- \mu_2$"].
\rm
Consider the parallel measurement
${\mathsf M}_{C_0(({\mathbb R} \times {\mathbb R}_+) \times ({\mathbb R} \times {\mathbb R}_+))}$
$({\mathsf O}_N^n \otimes {\mathsf O}_N^m= ({\mathbb R}^n \times {\mathbb R}^m \ , {\mathcal B}_{\mathbb R}^n \bigstimes {\mathcal B}_{\mathbb R}^m, {{{N}}^n}
\otimes  {{{N}}^m}) ,$
$S_{[(\mu_1, \sigma_1, \mu_2 , \sigma_2)]})$
in $C_0(({\mathbb R} \times {\mathbb R}_+) \times ({\mathbb R} \times {\mathbb R}_+))$.
\rm

Assume that $\sigma_1$ and $\sigma_2$ are fixed and known. Thus, this parallel 
measurement is represented by
${\mathsf M}_{C_0({\mathbb R} \times{\mathbb R} )}$
$({\mathsf O}_{N_{\sigma_1}}^n \otimes {\mathsf O}_{N_{\sigma_1}}^m= ({\mathbb R}^n \times {\mathbb R}^m \ , {\mathcal B}_{\mathbb R}^n \bigstimes {\mathcal B}_{\mathbb R}^m, {{{N_{\sigma_1}}}^n}
\otimes  {{{N_{\sigma_2}}}^m}) ,$
$S_{[(\mu_1,  \mu_2 )]})$
in $C_0({\mathbb R} \times {\mathbb R} )$. Here, recall the (\ref{eq2}),
i.e.,
\par
\noindent
\begin{align}
&
[{{{N_\sigma}}}({\Xi})] ({} {}{\mu} {}) 
=
\frac{1}{{\sqrt{2 \pi }\sigma{}}}
\int_{{\Xi}} \exp[{}- \frac{({}{}{x} - {}{\mu}  {})^2 }{2 \sigma^2}    {}] d {}{x}
\quad
({}\forall  {\Xi} \in {\cal B}_{{\mathbb R}{}}\mbox{(=Borel field in ${\mathbb R}$))},
\quad
\forall  \mu   \in {\mathbb R}).
\label{eq64}
\end{align}
Therefore, we have the state space
$\Omega ={\mathbb R}^2
=
\{
\omega=(\mu_1, \mu_2) \;:\; \mu_1,\mu_2
 \in {\mathbb R} 
\}$.
Put $\Theta={\mathbb R}$ with the distance $d_\Theta^{(1)} ( \theta_1, \theta_2 )= |\theta_1-\theta_2|$
and consider the quantity $\pi:{\mathbb R}^2 \to
{\mathbb R}$ by
\begin{align}
\pi (\mu_1, \mu_2)= \mu_1-\mu_2
\label{eq65}
\end{align}
The estimator $E: \widehat{X}(=X \times Y =
{{\mathbb R}^n \times {\mathbb R}^m}) 
\to \Theta(={\mathbb R})$ is defined by
\begin{align}
E(x_1, \ldots, x_n,y_1, \ldots, y_m)
=
\frac{\sum_{k=1}^n x_k}{n}
-
\frac{\sum_{k=1}^m y_k}{m}
\label{eq66}
\end{align}

For any
$ \omega=(\mu_1, \mu_2 )  ({}\in \Omega=
{\mathbb  R} \times {\mathbb R} )$,
define
the positive number
$\eta^\alpha_{\omega}$
$({}> 0)$
such that:
\begin{align}
\eta^\alpha_{\omega}
=
\inf
\{
\eta > 0:
[F ({}E^{-1} ({}
{{\rm Ball}^C_{d_\Theta^{(1)}}}(\pi(\omega) ; \eta{}))](\omega )
\ge \alpha
\}
\nonumber
\end{align}
where
${{\rm Ball}^C_{d_\Theta^{(1)} }}(\pi(\omega) ; \eta)$
$=
(-\infty, \mu_1 - \mu_2  - \eta] 
\cup
[ \mu_1 - \mu_2 + \eta , \infty)$.
Define the null hypothesis
$H_N$
$(\subseteq \Theta = {\mathbb R})$ such that
$$
H_N =\{ \theta_0 \}
$$

Now let us calculate the $\eta^\alpha_{\omega}$ as follows:
\begin{align}
&
E^{-1}({{\rm Ball}^C_{d_\Theta^{(1)} }}(\pi(\omega) ; \eta ))
=
E^{-1}(
(-\infty, \mu_1 - \mu_2  - \eta] 
\cup
[ \mu_1 - \mu_2 + \eta , \infty)
%
)
\nonumber
\\
=
&
\{
(x_1, \ldots , x_n, y_1, \ldots, y_m )
\in {\mathbb R}^n \times {\mathbb R}^m
\;:
\;
|
\frac{\sum_{k=1}^n x_k}{n}
-
\frac{\sum_{k=1}^m y_k}{m}
-(\mu_1 - \mu_2)|
\ge   \eta 
\}
\nonumber
\\
=
&
\{
(x_1, \ldots , x_n, y_1, \ldots, y_m )
\in {\mathbb R}^n \times {\mathbb R}^m
\;:
\;
|
\frac{\sum_{k=1}^n (x_k - \mu_1)}{n}
-
\frac{\sum_{k=1}^m (y_k- \mu_2)}{m}
|
\ge  \eta 
\}
\label{eq67}
\end{align}

Thus,
\begin{align}
&
[
({{{N_{\sigma_1}}}}^n
\otimes
{{{N_{\sigma_2}}}}^m
)
(E^{-1}({{\rm Ball}^C_{d_\Theta^{(1)} }}(\pi(\omega) ; \eta ))]
({}\omega{})
\nonumber
\\
=
&
\frac{1}{({{\sqrt{2 \pi }\sigma_1{}}})^n({{\sqrt{2 \pi }\sigma_2{}}})^m}
\nonumber
\\
&
\bigtimes 
\!\!\!\!\!\!
\underset{{
|
\frac{\sum_{k=1}^n( x_k - \mu_1)}{n}
-
\frac{\sum_{k=1}^m (y_k- \mu_2)}{m}
|
\ge  \eta 
%
%
%
}}{\int \cdots \int}
\exp[
{}- \frac{\sum_{k=1}^n ({}{}{x_k} - {}{\mu_1}  {})^2 
}
{2 \sigma_1^2}
{}- \frac{\sum_{k=1}^m ({}{}{y_k} - {}{\mu_2}  {})^2 
}
{2 \sigma_2^2}
] d {}{x_1} d {}{x_2}\cdots dx_nd {}{y_1} d {}{y_2}\cdots dy_m
\nonumber
\\
=
&
\frac{1}{({{\sqrt{2 \pi }\sigma_1{}}})^n({{\sqrt{2 \pi }\sigma_2{}}})^m}
\underset{{
|
\frac{\sum_{k=1}^n x_k }{n}
-
\frac{\sum_{k=1}^m y_k}{m}
|
\ge  \eta 
%
%
%
}}{\int \cdots \int}
\exp[
- \frac{
\sum_{k=1}^n {x_k}^2 
}
{2 \sigma_1^2}
- \frac{
\sum_{k=1}^m {y_k}^2 
}
{2 \sigma_2^2}
] d {}{x_1} d {}{x_2}\cdots dx_nd {}{y_1} d {}{y_2}\cdots dy_m
\nonumber
\\
=
&
1-
\frac{1}{{\sqrt{2 \pi }(\frac{\sigma_1^2}{n}+\frac{\sigma_2^2}{m})^{1/2}{}}}
\int_{{- \eta}}^{\eta} \exp[{}- \frac{{x}^2 }{2 (\frac{\sigma_1^2}{n}+\frac{\sigma_2^2}{m})}] d {x}
\label{eq68}
\end{align}
Using the $z(\alpha/2)$ in (\ref{eq33}),
we get that
\begin{align}
\eta^\alpha_{\omega} = 
(\frac{\sigma_1^2}{n}+\frac{\sigma_2^2}{m})^{1/2}
z(\frac{\alpha}{2})
\label{eq69}
\end{align}
Therefore,
we get ${\widehat R}_{\widehat{x}}^{\alpha}$
(
the $({}\alpha{})$-rejection region
of
$H_N =\{\theta_0\}( \subseteq \Theta)$
)
as follows:
\begin{align}
{\widehat R}_{H_N}^{\alpha}
&
=
\bigcap_{\omega =(\mu_1, \mu_2 ) \in  \Omega (={\mathbb R}^2 ) \mbox{ \footnotesize such that }
\pi(\omega)= \mu_1-\mu_2 \in {H_N}(=\{\theta_0 \} )}
\{
E(\widehat{x})
(\in
\Theta)
:
d_\Theta^{(1)}  ({}E(\widehat{x}),
\pi(\omega))
\ge
\eta^\alpha_{\omega }
\}
\nonumber
\\
&
=
\{ \overline{\mu}(x)-\overline{\mu}(y) \in \Theta (={\mathbb R})
\;:\;
|
\overline{\mu}(x)-\overline{\mu}(y)
-\theta_0|
\ge 
(\frac{\sigma_1^2}{n}+\frac{\sigma_2^2}{m})^{1/2}
 z(\frac{\alpha}{2})
 \}
\label{eq70}
\end{align}
where
$$
\overline{\mu}(x)=\frac{\sum_{k=1}^n x_k }{n},
\quad
\overline{\mu}(y)
=
\frac{\sum_{k=1}^m y_k}{m}
$$
%
\rm
\par
\noindent
\bf
Remark 3
\rm
[The case that $H_N = (- \infty , \theta_0]$].
If the null hypothesis
$H_N$ is assumed as follows:
$$
H_N = (- \infty , \theta_0],
$$
it suffices to define the semi-distance
\begin{align}
d_\Theta^{(1)}  (\theta_1, \theta_2)
=
\cases
|\theta_1-\theta_2|
\quad
&
(
\forall \theta_1, \theta_2 \in \Theta={\mathbb R}
\mbox{ such that }
\theta_0 \le \theta_1, \theta_2
)
\\
\max \{ \theta_1, \theta_2 \}
- \theta_0
\quad
&
(
\forall \theta_1, \theta_2 \in \Theta={\mathbb R}
\mbox{ such that }
\min \{ \theta_1, \theta_2 \}
\le
\theta_0 \le \max \{ \theta_1, \theta_2 \}
)
\\
0
&
(
\forall \theta_1, \theta_2 \in \Theta={\mathbb R}
\mbox{ such that }
\theta_1, \theta_2
\le \theta_0
)
\endcases
\label{eq71}
\end{align}
Then, we can easily see that
\begin{align}
{\widehat R}_{H_N}^{\alpha}
&
=
\bigcap_{\omega =(\mu_1, \mu_2 ) \in  \Omega (={\mathbb R}^2 ) \mbox{ \footnotesize such that }
\pi(\omega)= \mu_1-\mu_2 \in {H_N}(=(- \infty, \theta_0] )}
\{
E(\widehat{x})
(\in
\Theta)
:
d_\Theta^{(1)}  ({}E(\widehat{x}),
\pi(\omega))
\ge
\eta^\alpha_{\omega }
\}
\nonumber
\\
&
=
\{ \overline{\mu}(x)-\overline{\mu}(y) \in {\mathbb R} 
\;:\;
\overline{\mu}(x)-\overline{\mu}(y)
-\theta_0
\ge 
(\frac{\sigma_1^2}{n}+\frac{\sigma_2^2}{m})^{1/2}
 z({\alpha}{})
 \}
\label{eq72} 
\end{align}

\subsection{The ratio of the population variances
}
\par
\noindent
\rm
\par
\noindent
\bf
Example 9
\rm
\rm
[Rejection region in the case that "$\pi(\sigma_1, \sigma_2)=\sigma_1
/\mu_2$"].
\rm
Consider the parallel measurement
${\mathsf M}_{C_0(({\mathbb R} \times {\mathbb R}_+) \times ({\mathbb R} \times {\mathbb R}_+))}$
$({\mathsf O}_N^n \otimes {\mathsf O}_N^m= 
({\mathbb R}^n \times {\mathbb R}^m \ , {\mathcal B}_{\mathbb R}^n \bigstimes {\mathcal B}_{\mathbb R}^m, {{{N}}^n}
\otimes  {{{N}}^m}) ,$
$S_{[(\mu_1, \sigma_1, \mu_2 , \sigma_2)]})$
in $C_0(({\mathbb R} \times {\mathbb R}_+) \times ({\mathbb R} \times {\mathbb R}_+))$.
\rm

Put $\Theta={\mathbb R}_+$ with the distance $d_\Theta^{(2)} ( \theta_1, \theta_2 )= |\log{\theta_1}-\log{\theta_2}|=|\log{\frac{\theta_1}{\theta_2}}|$
and consider the quantity $\pi:\Omega = ({\mathbb R}\times {\mathbb R}_+) \times
({\mathbb R}\times {\mathbb R}_+)
\to
\Theta= {\mathbb R}_+$ by
\begin{align}
\pi ((\mu_1,\sigma_1), (\mu_2, \sigma_2))= \sigma_1/\sigma_2
\label{eq73}
\end{align}
The estimator $E: \widehat{X}(=X \times Y =
{{\mathbb R}^n \times {\mathbb R}^m}) 
\to \Theta(={\mathbb R}_+)$ is defined by
\begin{align}
E(x_1, \ldots, x_n,y_1, \ldots, y_m)
=
\frac{
\overline{\sigma}'_1(x)
}
{\overline{\sigma}'_2(y)
}
\qquad
\mbox{(Recall ((\ref{eq16}))}
%
\label{eq74}
\end{align}

For any
$\omega =((\mu_1,\sigma_1), (\mu_2, \sigma_2))
\in
\Omega = ({\mathbb R}\times {\mathbb R}_+) \times
({\mathbb R}\times {\mathbb R}_+)$,
define
the positive number
$\eta^\alpha_{\omega}$
$({}> 0)$
such that:
\begin{align}
\eta^\alpha_{\omega}
=
\inf
\{
\eta > 0:
[F ({}E^{-1} ({}
{{\rm Ball}^C_{d_\Theta^{(2)} }}(\pi(\omega) ; \eta{}))](\omega )
\ge \alpha
\}
\nonumber
\end{align}
where
${{\rm Ball}^C_{d_\Theta^{(2)} }}(\pi(\omega) ; \eta)$
$=
(0, (\sigma_1/\sigma_2)e^{-\eta} ] 
\cup
[ (\sigma_1/\sigma_2)e^{\eta} , \infty)$.
Define the null hypothesis
$H_N$
$(\subseteq \Theta = {\mathbb R}_+)$ such that
$$
H_N =\{ r_0 \}
$$

Now let us calculate the $\eta^\alpha_{\omega}$ as follows:
\begin{align}
&
E^{-1}({{\rm Ball}^C_{d_\Theta^{(2)} }}(\pi(\omega) ; \eta ))
=
E^{-1}(
(0, (\sigma_1/\sigma_2)e^{-\eta} ] 
\cup
[ (\sigma_1/\sigma_2)e^{\eta} , \infty)
)
\nonumber
\\
=
&
\{
(x_1, \ldots , x_n, y_1, \ldots, y_m )
\in {\mathbb R}^n \times {\mathbb R}^m
\;:
\;
\frac{
\overline{\sigma}'_1(x)/ \sigma_1
}
{\overline{\sigma}'_2(y)/ \sigma_2
}
\le e^{-\eta}
\mbox{ or }
\frac{
\overline{\sigma}'_1(x)/ \sigma_1
}
{\overline{\sigma}'_2(y)/ \sigma_2
}
\ge e^{\eta}
\}
\label{eq75}
\end{align}

Thus,
\begin{align}
&
1-
[
({{{N_{\sigma_1}}}}^n
\otimes
{{{N_{\sigma_2}}}}^m
)
(E^{-1}({{\rm Ball}^C_{d_\Theta^{(2)} }}(\pi(\omega) ; \eta ))]
({}\omega{})
\nonumber
\\
=
&
\frac{1}{({{\sqrt{2 \pi }\sigma_1{}}})^n({{\sqrt{2 \pi }\sigma_2{}}})^m}
\nonumber
\\
&
\bigtimes 
\!\!\!\!\!\!
\underset{{
e^{- \eta}
\le
\frac{
\overline{\sigma}'_1(x)/ \sigma_1
}
{\overline{\sigma}'_2(y)/ \sigma_2
}
\le
e^{\eta}
}}{\int \cdots \int}
\exp[
{}- \frac{\sum_{k=1}^n ({}{}{x_k} - {}{\mu_1}  {})^2 
}
{2 \sigma_1^2}
{}- \frac{\sum_{k=1}^m ({}{}{y_k} - {}{\mu_2}  {})^2 
}
{2 \sigma_2^2}
] d {}{x_1} d {}{x_2}\cdots dx_nd {}{y_1} d {}{y_2}\cdots dy_m
\nonumber
\\
=
&
\frac{1}{({{\sqrt{2 \pi }{}}})^n({{\sqrt{2 \pi }{}}})^m}
\underset{{
e^{- \eta}
\le
\frac{
\overline{\sigma}'_1(x)
}
{\overline{\sigma}'_2(y)
}
\le
e^{\eta}
}}{\int \cdots \int}
\exp[
- \frac{
\sum_{k=1}^n {x_k}^2 
}
{2 }
- \frac{
\sum_{k=1}^m {y_k}^2 
}
{2 }
] d {}{x_1} d {}{x_2}\cdots dx_nd {}{y_1} d {}{y_2}\cdots dy_m
\nonumber
\\
=
&
\int_{e^{-2\eta}}^{e^{2\eta}}
p^F_{n-1,m-1}(x) dx
\label{eq76}
\end{align}
where
$p^{F}_{n-1,m-1}({ x} )$ is the $F$-distribution with $(n-1,m-1)$ degrees of freedom.
Define the positive $\eta^\alpha_{\omega}$
such that
$$
1- \alpha
=
\int_{e^{-2\eta^\alpha_{\omega}}}^{e^{2\eta^\alpha_{\omega}}}
p^F_{n-1,m-1}(x) dx
$$
Since
$\eta^\alpha_{\omega}$ does not depemd on
$\omega$, we can put
$\eta^\alpha_{n}$
$= \eta^\alpha_{\omega}$.
Therefore,
we get ${\widehat R}_{H_N}^{\alpha}$(
\it
the $({}\alpha{})$-rejection region
of
$H_N(=\{ r \})$
\rm
)
as follows:
\begin{align}
{\widehat R}_{H_N}^{\alpha}
&
=
\bigcap_{\omega =(\mu_1, \sigma_1,\mu_2, \sigma_2 ) \in  \Omega
\mbox{ \footnotesize such that }
\pi(\omega)= \sigma_1/\sigma_2 \in {H_N}(=\{ r_0 \} )}
\{
E(\widehat{x})
(\in
\Theta)
:
d_\Theta^{(2)}  ({}E(\widehat{x}),
\pi(\omega))
\ge
\eta^\alpha_{n }
\}
\nonumber
\\
&
=
\bigcap_{
{\frac{\sigma_1}{\sigma_2}=r_0}
}
\{
\frac{
\overline{\sigma}'_1(x)
}
{\overline{\sigma}'_2(y)
}
\in \Theta = {\mathbb R}_+
\;:
\;
\frac{
\overline{\sigma}'_1(x)/ \sigma_1
}
{\overline{\sigma}'_2(y)/ \sigma_2
}
\le e^{-\eta^\alpha_{n }}
\mbox{ or }
\frac{
\overline{\sigma}'_1(x)/ \sigma_1
}
{\overline{\sigma}'_2(y)/ \sigma_2
}
\ge e^{\eta^\alpha_{n }}
\}
\nonumber
\\
&
=
\{
\frac{
\overline{\sigma}'_1(x)
}
{\overline{\sigma}'_2(y)
}
\in \Theta = {\mathbb R}_+
\;:
\;
\frac{
\overline{\sigma}'_1(x)
}
{\overline{\sigma}'_2(y)
}
\le r_0 e^{-\eta^\alpha_{n }}
\mbox{ or }
\frac{
\overline{\sigma}'_1(x)
}
{\overline{\sigma}'_2(y)
}
\ge r_0 e^{\eta^\alpha_{n }}
\}
\label{eq77}
\end{align}
%
\rm
\par
\noindent
\bf
Remark 4
\rm
[The case that $H_N = (0 , r_0] \subseteq \Theta = {\mathbb R}_+$].
If the null hypothesis
$H_N$ is assumed as follows:
$$
H_N = (0 , r_0],
$$
it suffices to define the semi-distance
\begin{align}
d_\Theta^{(2)}  (\theta_1, \theta_2)
=
\cases
|\log(\theta_1/\theta_2)|
\quad
&
(
\forall \theta_1, \theta_2 \in \Theta={\mathbb R}
\mbox{ such that }
r_0 \le \theta_1, \theta_2
)
\\
|\log(\max \{ \theta_1, \theta_2 \}/r_0)
\quad
&
(
\forall \theta_1, \theta_2 \in \Theta={\mathbb R}
\mbox{ such that }
\min \{ \theta_1, \theta_2 \}
\le
r_0 \le \max \{ \theta_1, \theta_2 \}
)
\\
0
&
(
\forall \theta_1, \theta_2 \in \Theta={\mathbb R}
\mbox{ such that }
\theta_1, \theta_2
\le r_0
)
\endcases
\label{eq78}
\end{align}
Then, we can easily see that
\begin{align}
{\widehat R}_{H_N}^{\alpha}
&
=
\bigcap_{\omega =(\mu_1, \mu_2 ) \in  \Omega (={\mathbb R}^2 ) \mbox{ \footnotesize such that }
\pi(\omega)= \sigma_1 / \sigma_2 \in {H_N}(=(0, r_0] )}
\{
E(\widehat{x})
(\in
\Theta)
:
d_\Theta^{(2)}  ({}E(\widehat{x}),
\pi(\omega))
\ge
\eta^\alpha_{\omega}
\}
\nonumber
\\
&
=
\{
\frac{
\overline{\sigma}'_1(x)
}
{\overline{\sigma}'_2(y)
}
\in \Theta = {\mathbb R}_+
\;:
\;
\frac{
\overline{\sigma}'_1(x)
}
{\overline{\sigma}'_2(y)
}
\ge r_0 e^{(\eta^\alpha_{n })'}
\}
%
%
%
%
%
%
\nonumber
\\
&
=
[ r_0 e^{(\eta^\alpha_{n })'} , \infty )
\label{eq79} 
\end{align}
where the positive $(\eta^\alpha_{n})'$
such that
$$
\alpha
=
\int_{e^{2 (\eta^\alpha_{n})'}}^{\infty}
p^F_{n-1,m-1}(x) dx
$$

\subsection{The case that
$d^x_\Theta$ depends on $x$;
Student's t-distribution}
\par
\noindent
\rm
\par
The arguments in this section are continued from Example 2.
\par
\noindent
\bf
Example 10
\rm
[Student's t-distribution].
Consider the simultaneous measurement
${\mathsf M}_{C_0({\mathbb R} \times {\mathbb R}_+)}$
$({\mathsf O}_N^n = ({\mathbb R}^n, {\mathcal B}_{\mathbb R}^n, {{{N}}^n}) ,$
$S_{[(\mu, \sigma)]})$
in $C_0({\mathbb R} \times {\mathbb R}_+)$.
Thus,
we consider that
$\Omega = {\mathbb R} \times {\mathbb R}_+$,
$X={\mathbb R}^n$.
Put
$\Theta={\mathbb R}$ with the semi-distance
$d_\Theta^x (\forall x \in X)$
such that
\begin{align}
d_\Theta^x (\theta_1, \theta_2)
=
\frac{|\theta_1-\theta_2|}{{\overline{\sigma}'(x)}/\sqrt{n}}
\quad
\qquad
(\forall x \in X={\mathbb R}^n,
\forall \theta_1, \theta_2 \in \Theta={\mathbb R}
)
\label{eq80}
\end{align}
where ${\overline{\sigma}'(x)}=\sqrt{\frac{n}{n-1}}\overline{\sigma}(x)$.
The quantity $\pi:\Omega(={\mathbb R} \times {\mathbb R}_+)
\to
\Theta(={\mathbb R})$
is defined by
\begin{align}
\Omega(={\mathbb R} \times {\mathbb R}_+)
\ni \omega
=
(\mu, \sigma )
\mapsto \pi (\mu, \sigma )
=
\mu
\in
\Theta(={\mathbb R})
\label{eq81}
\end{align}
\rm
Also, define the estimator
$E:X(={\mathbb R}^n) \to \Theta(={\mathbb R})$
such that
\begin{align}
E(x)=E(x_1, x_2, \ldots , x_n )
=
\overline{\mu}(x)
=
\frac{x_1 + x_2 + \cdots + x_n}{n}
\label{eq82}
\end{align}
Define the null hypothesis $H_N$
$(\subseteq
\Theta=
{\mathbb  R}  )
)$
such that
\begin{align}
H_N=  \{\mu_0\}
\label{eq83}
\end{align}
Thus, for any
$ \omega=(\mu_0, \sigma )  ({}\in \Omega=
{\mathbb  R} \times {\mathbb R}_+ )$,
we see that
\begin{align}
&
[N^n(\{ x \in X \;:\; 
d^x_\Theta ( E(x) , \pi( \omega ) )
\ge \eta
\}
)](\omega )
\nonumber
\\
=&
[N^n(\{ x \in X \;:\; 
\frac{
|\overline{\mu}(x)- \mu_0 |}{
{{\overline{\sigma}'(x)}/\sqrt{n}}
}
\ge \eta
\}
)](\omega )
\nonumber
\\
=
&
\frac{1}{({{\sqrt{2 \pi }\sigma{}}})^n}
\underset{
\eta
\le
\frac{
|\overline{\mu}(x)- \mu_0 |}{
{{\overline{\sigma}'(x)}/\sqrt{n}}
}
}{\int \cdots \int}
\exp[{}- \frac{\sum_{k=1}^n ({}{}{x_k} - {}{\mu_0}  {})^2 
}
{2 \sigma^2}    {}] d {}{x_1} d {}{x_2}\cdots dx_n
\nonumber
\\
=
&
\frac{1}{({{\sqrt{2 \pi }{}}})^n}
\underset{
\eta
\le
\frac{
|\overline{\mu}(x) |}{
{{\overline{\sigma}'(x)}/\sqrt{n}}
}
}
{\int \cdots \int}
\exp[{}- \frac{\sum_{k=1}^n ({}{}{x_k}  {}  {})^2 
}
{2 }    {}] d {}{x_1} d {}{x_2}\cdots dx_n
\nonumber
\\
\noindent
=
&
1-
\int_{-\eta}^{\eta
}
p^t_{n-1}(x)
dx
\label{eq84}
\end{align}
where
$p^t_{n-1}$
is
the t-distribution with $n-1$ degrees of freedom.
Solving the equation
$
1-\alpha
=
\int_{-\eta^\alpha_{\omega}}^{\eta^\alpha_{\omega}
}
p^t_{n-1}(x)
dx
$,
we get $\eta^\alpha_{\omega}$
$=t(\alpha/2)$.

%
%
%
%
%
%
%

Therefore,
we get ${\widehat R}_{H_N}^{\alpha}$(
\it
the $({}\alpha{})$-rejection region
of
$H_N(=\{\mu_0\})$
\rm
)
as follows:
\begin{align}
{\widehat R}_{H_N}^{\alpha}
&
=
\bigcap_{\omega =(\mu, \sigma ) \in  \Omega (={\mathbb R} \times {\mathbb R}_+) \mbox{ \footnotesize such that }
\pi(\omega)= \mu \in {H_N}(=\{\mu_0\})}
\{
E(x)
(\in
\Theta)
:
\;\;
d^x_\Theta ({}E(x),
\pi(\omega )
)
\ge
\eta^\alpha_{\omega }
\}
\nonumber
\\
&
=
\{\overline{\mu}(x) \in \Theta(={\mathbb R})
\;:\;
\frac{
|\overline{\mu}(x)- \mu_0 |}{
{{\overline{\sigma}'(x)}/\sqrt{n}}
}
\ge
t(\alpha/2)
\}
\nonumber
\\
&
=
\{\overline{\mu}(x) \in \Theta(={\mathbb R})
\;:\;
\mu_0 \le 
\overline{\mu}(x)
-
\frac{{\overline{\sigma}'(x)}}{\sqrt{n}}
t(\alpha/2)
\mbox{ or }
\overline{\mu}(x)
+
\frac{{\overline{\sigma}'(x)}}{\sqrt{n}}
t(\alpha/2)
\le \mu_0
\}
\label{eq85} 
\end{align}
\rm
\par
\noindent
\bf
Remark 5
\rm
[The case that $H_N = (- \infty , \mu_0]$].
If the null hypothesis
$H_N$ is assumed as follows:
$$
H_N = (- \infty , \mu_0],
$$
it suffices to define the semi-distance
\begin{align}
d_\Theta^x (\theta_1, \theta_2)
=
\cases
\frac{|\theta_1-\theta_2|}{{\overline{\sigma}'(x)}/\sqrt{n}}
\quad
&
(
\forall \theta_1, \theta_2 \in \Theta={\mathbb R}
\mbox{ such that }
\mu_0 \le \theta_1, \theta_2
)
\\
\frac{\max \{ \theta_1, \theta_2 \}-\mu_0}{{\overline{\sigma}'(x)}/\sqrt{n}}
\quad
&
(
\forall \theta_1, \theta_2 \in \Theta={\mathbb R}
\mbox{ such that }
\min \{ \theta_1, \theta_2 \}
\le
\mu_0 \le \max \{ \theta_1, \theta_2 \}
)
\\
0
&
(
\forall \theta_1, \theta_2 \in \Theta={\mathbb R}
\mbox{ such that }
\theta_1, \theta_2
\le \mu_0
)
\endcases
\label{eq86}
\end{align}
for any
$x \in X={\mathbb R}^n$.
Then, we can easily see that
\begin{align}
{\widehat R}_{H_N}^{\alpha}
&
=
\bigcap_{\omega =(\mu, \sigma ) \in  \Omega (={\mathbb R} \times {\mathbb R}_+) \mbox{ \footnotesize such that }
\pi(\omega)= \mu \in {H_N}(=(- \infty , \mu_0])}
\{
E(x)
(\in
\Theta)
:
\;\;
d^x_\Theta ({}E(x),
\pi(\omega )
)
\ge
\eta^\alpha_{\omega }
\}
\nonumber
\\
&
=
\{\overline{\mu}(x) \in \Theta(={\mathbb R})
\;:\;
\mu_0 \le 
\overline{\mu}(x)
-
\frac{{\overline{\sigma}'(x)}}{\sqrt{n}}
t(\alpha)
\}
\label{eq87} 
\end{align}

\section{Conclusions
}
\par
\noindent
\par
It is sure that
statistics and (classical) quantum language
are similar. however,
quantum language has the firm structure (\ref{eq1}),
i.e.,
\begin{align}
\underset{\mbox{(=MT(measurement theory))}}{\fbox{Quantum language}}
=
\underset{\mbox{(measurement)}}{\fbox{Axiom 1}}
+
\underset{\mbox{(causality)}}{\fbox{Axiom 2}}
+
\underset{\mbox{(how to use Axioms)}}{\fbox{linguistic interpretation}}
\label{eq88}
\end{align}
Hence,
as seen in Theorems 1-4 of this paper,
every argument cannot but become clear
in quantum language.

\par
\noindent
\par
Particularly,
the following two statistical hypothesis tests
(J$_1$) and (J$_2$), that is, 
\begin{itemize}
\item[(J$_1$)] Theorem 2 (Likelihood ratio test)
\\
key-words: Estimator $E:X \to \Omega$, Quantity $\pi:\Omega \to \Theta$, \fbox{Likelihood function $L_\theta(\omega)$ in (\ref{eq17})}
\item[(J$_2$)] Theorem 4 (Reverse confidence interval method)
\\
key-words: Estimator $E:X \to \Theta$, Quantity $\pi:\Omega \to \Theta$, \fbox{Semi-distance $d^x_\Theta$ on $\Theta$}.
\end{itemize}
should be compared and examined. 

\par
For example, we remark that the difference between "one sided test" and "two sided test" is due to the difference of the semi-distances. 
And further, we see the peculiarity of the student's $t$-distribution in Example 10,
however, we have no firm answer to the following question:
\begin{itemize}
\item[(K)] Can Example 10 (Student's t-distribution)
be naturally understood in Theorem 2 (Likelihood ratio test)?
\end{itemize}
Although Theorem 2 (Likelihood ratio test) is orthodox, it is not handy.
On the other hand, we believe that Theorems 4 (Reverse confidence interval method)
may be usual,
though it is not presented as a general theorem in the elementary books of statistics.

\vskip0.3cm

Since quantum language is suited for theoretical arguments, we believe, from the theoretical point of view, that our results
(i.e.,
Theorems 1-4) are final in classical systems.
We hope that our assertions will be examined from various points of view.



\rm
\par
\renewcommand{\refname}{
\large 
References}
{
\small

\normalsize
}


\begin{thebibliography}{9}
\rm

\rm
\bibitem{Davi} E. B. Davies,
\newblock {\em Quantum Theory of Open Systems,}
\newblock {Academic Press,}
\newblock {1976}
\bibitem{Ishi2}
S. Ishikawa,
{\em A Quantum Mechanical Approach to Fuzzy Theory,}
{{\rm Fuzzy Sets and Systems}}, 
{Vol. 90}, No. 3, {} 277-306,
1997
\\
\href{http://dx.doi.org/10.1016/S0165-0114(96)00114-5}{doi: 10.1016/S0165-0114(96)00114-5}
\bibitem{Ishi3}
S. Ishikawa,
{\it Statistics in measurements}, 
Fuzzy sets and systems, 
{Vol. 116}, No. 2, 141-154,
2000
\\
\href{http://dx.doi.org/10.1016/S0165-0114(98)00280-2}{doi:10.1016/S0165-0114(98)00280-2}
\bibitem{Ishi4}
S. Ishikawa,
{\em Mathematical Foundations of Measurement Theory,}
Keio University Press Inc. 335pages,
2006.
\\
(\url{http://www.keio-up.co.jp/kup/mfomt/})
\bibitem{Ishi5} {S. Ishikawa,}
\newblock {\em A New Interpretation of Quantum Mechanics,}
\newblock {\rm Journal of quantum information science},
{Vol. 1}, No. 2, {}35-42,
2011
\\
\href{http://dx.doi.org/10.4236/jqis.2011.12005}{doi: 10.4236/jqis.2011.12005}
\\
(\url{http://www.scirp.org/journal/PaperInformation.aspx?paperID=7610})
\bibitem{Ishi6} {S. Ishikawa,}
\newblock {\em Quantum Mechanics and the Philosophy of Language:
Reconsideration of traditional
philosophies,}
\newblock {\rm Journal of quantum information science},
{Vol. 2}, No. 1, {}2-9,
2012
\\
\href{http://dx.doi.org/ 10.4236/jqis.2012.21002}{doi: 10.4236/jqis.2012.21002}
\\
(\url{
http://www.scirp.org/journal/PaperInformation.aspx?paperID=18194
})
\bibitem{Ishi7} {S. Ishikawa,}
\newblock {\em A Measurement Theoretical
Foundation of Statistics,}
\newblock {\rm Journal of Applied Mathematics},
{Vol. 3}, No. 3, {} 283-292,
2012
\\
\href{http://dx.doi.org/10.4236/am.2012.33044}{doi: 10.4236/am.2012.33044}
\\
(\url{http://www.scirp.org/journal/PaperInformation.aspx?paperID=18109&})
\bibitem{Ishi8} S. Ishikawa,
\newblock {\em What is statistics? The answer by quantum language,}
\href{http://arxiv.org/abs/1207.0407}{arXiv:1312.6757 [math.ST]}
\newblock {2012.}
\\
(\url{http://arxiv.org/abs/1207.0407})
\bibitem{Ishi9} S. Ishikawa, K. Kikuchi,
\newblock {\em The Confidence Interval Methods in Quantum Language,}
\href{http://arxiv.org/abs/1312.6757}{arXiv:1312.6757 [math.ST]}
\newblock {2013.}
\\
(\url{http://arxiv.org/abs/1312.6757})
\bibitem{Neum} {J. von Neumann,}
\newblock {\em Mathematical Foundations of Quantum Mechanics,}
\newblock {\rm Springer Verlag, Berlin,}
\newblock {1932}
\bibitem{Saka}
S. Sakai,
{\it $C^*$-algebras and $W^*$-algebras}, 
Ergebnisse der Mathematik und ihrer Grenzgebiete (Band 60), 
Springer-Verlag, (1971)
\bibitem{Yosi}
K. Yosida,
{\em Functional Analysis,}
Springer-Verlag, 6th edition, 1980.
\end{thebibliography}
\end{document}